\documentclass[12pt]{article}

\usepackage{amssymb,amsmath,amsthm}

\newcommand{\besovspq}{\ensuremath{B ^s _{p,q}}}
\newcommand{\bord}[1]{\ensuremath{\partial#1}}
\newcommand{\caract}[1]{\ensuremath{{\bf 1}_{#1}}}
\newcommand{\cinf}{\ensuremath{C^\infty}}
\newcommand{\cinfc}{\ensuremath{C^\infty_0}}
\newcommand{\conv}[2]{\ensuremath{#1\ast#2}}
\newcommand{\Cr}{\ensuremath{C^{r}}}
\newcommand{\cronde}{\ensuremath{\mathcal C}}
\newcommand{\Crp}{\ensuremath{C^{r'}}}

\newcommand{\cuet}{\ensuremath{C^{1}_{\etoile}}}
\newcommand{\cun}{\ensuremath{C^{1}}}
\newcommand{\Cunetoile}{\cuet}
\newcommand{\cunplusr}{\ensuremath{C^{1+r}}}

\newcommand{\deln}{\ensuremath{\Delta_{n}}}
\newcommand{\delp}{\ensuremath{\Delta_{p}}}
\newcommand{\delq}{\ensuremath{\Delta_{q}}}
\newcommand{\deltap}{\delp}
\newcommand{\deltaq}{\delq}
\newcommand{\demonstration}{\begin{proof}}
\newcommand{\deqa}{\begin{eqnarray*}}

\newcommand{\distance}{\ensuremath{\mbox{\rm d}}}
\newcommand{\divergence}{\ensuremath{\mbox{\rm div}\,}}

\newcommand{\domec}{\ensuremath{\Omega}}
\newcommand{\dsdn}[1]{\ensuremath{\frac{\partial#1}{\partial n}}}
\newcommand{\egdef}{\ensuremath{\stackrel{\mathrm{def}}{=}}}
\newcommand{\egnot}{\ensuremath{\stackrel{\mathrm{not}}{=}}}
\newcommand{\etoile}{\ensuremath{\star}}
\newcommand{\fcaract}[1]{\caract{#1}}
\newcommand{\feqa}{\end{eqnarray*}}
\newcommand{\ferm}{\fermeture}
\newcommand{\fermeture}[1]{\ensuremath{\overline{#1}}}
\newcommand{\fermeturedans}{\ensuremath{\subset\subset}}

\newcommand{\findemonstration}{\end{proof}}
\newcommand{\fonctioncaracteristique}[1]{\caract{#1}}

\newcommand{\lap}{\ensuremath{\Delta}}
\newcommand{\Linf}{\ensuremath{L^\infty}}
\newcommand{\lip}{\Lip}
\newcommand{\Lip}{\ensuremath{\mathrm{Lip}}}
\newcommand{\Lp}{\ensuremath{L^{p}}}

\newcommand{\naturels}{\nn}
\newcommand{\nL}[1]{\ensuremath{\|#1\|_{L^{1}}}}
\newcommand{\nLinf}[1]{\ensuremath{\|#1\|_{\Linf}}}
\newcommand{\nlip}[1]{\ensuremath{\norme{#1}_\lip}}
\newcommand{\nLp}[1]{\ensuremath{\|#1\|_{\Lp}}}
\newcommand{\nn}{\ensuremath{{\mathbb N}}}
\newcommand{\norme}[1]{\ensuremath{\|#1\|}}
\newcommand{\normeCr}{\nr}
\newcommand{\normeCrmu}{\nrmi}
\newcommand{\normeCrpmu}{\nrpmi}
\newcommand{\normeCs}{\ns}
\newcommand{\normeCsmu}{\nsmi}
\newcommand{\normeLip}{\nlip}
\newcommand{\normeLinf}{\nLinf}
\newcommand{\normeLun}{\nL}
\newcommand{\normeun}[1]{\ensuremath{\|#1\|_{1}}}
\newcommand{\normezer}[1]{\ensuremath{\|#1\|_{0}}}
\newcommand{\nr}[1]{\ensuremath{\|#1\|_{r}}}
\newcommand{\nrmi}[1]{\ensuremath{\|#1\|_{r-1}}}

\newcommand{\nrpmi}[1]{\ensuremath{\|#1\|_{r'-1}}}
\newcommand{\ns}[1]{\ensuremath{\|#1\|_{s}}}
\newcommand{\nsmi}[1]{\ensuremath{\|#1\|_{s-1}}}

\newcommand{\om}{\ensuremath{\omega}}
\newcommand{\omegazero}{\omzer}
\newcommand{\omzer}{\ensuremath{\om_{0}}}
\newcommand{\paraproduit}[2]{\ensuremath{T_{#1}#2}}
\newcommand{\pgrad}{\ensuremath{\cdot\nabla}}
\newcommand{\poche}{\ensuremath{P}}
\newcommand{\produitcartesien}{\times}
\newcommand{\produitvectoriel}{\times}

\newcommand{\rdeux}{\ensuremath{{\mathbb R}^2}}
\newcommand{\reste}[2]{\ensuremath{R(#1,\,#2)}}
\newcommand{\restr}[2]{\ensuremath{\left. #1 \right|_{#2}}}
\newcommand{\rond}{\ensuremath{\circ}}

\newcommand{\rotationnel}{\ensuremath{\mbox{\rm rot}\,}}
\newcommand{\rtrois}{\ensuremath{{\mathbb R}^3}}
\newcommand{\rr}{\ensuremath{{\mathbb R}}}
\newcommand{\schw}{\ensuremath{\mathcal{S}}}

\newcommand{\support}[1]{\ensuremath{\mathrm{supp}\ #1}}
\newcommand{\teste}[2]{\ensuremath{\left\langle #1,\,#2\right\rangle}}
\newcommand{\uronde}{\ensuremath{\mathcal U}}
\newcommand{\urondep}{\ensuremath{\mathcal U'}}
\newtheorem{thm}{Theorem}
\newtheorem{definition}{Definition}
\newcommand{\trace}{\ensuremath{\mathrm{tr}}}

\newcommand{\vinitial}{\ensuremath{v_0}}

\newtheorem{prop}{Proposition}
\newtheorem{lem}{Lemma}

\begin{document}
\mbox{}
\vspace{2in}
\begin{center}
\textbf{ON 3-D VORTEX PATCHES IN BOUNDED \vspace{2\baselineskip}DOMAINS\footnotemark}\\
\textbf{Alexandre \vspace{\baselineskip}Dutrifoy}\\
Laboratoire Jacques-Louis Lions\\
Universit{\'e} Pierre et Marie Curie\\
Bo{\^{\i}}te courrier 187\\
75252 Paris Cedex 05\\
France
\end{center}
\vspace{4\baselineskip}

\section{Introduction}

\subsection{Abstract}

This article concerns the equations of motion of perfect incompressible fluids in a smooth,
bounded, simply connected domain of \rtrois. So we study the Euler system
\begin{equation}
\left\{\begin{array}{l}
\partial_t v + v \pgrad v = -\nabla p, \\
\divergence v = 0,
\end{array}\right.\label{E}
\end{equation}
where $v$ is the velocity field and $p$ is the pressure, along with an initial datum,
\begin{equation}
\restr{v }{t=0 } = \vinitial ,\label{I}
\end{equation}
and a condition at the boundary of the domain $\Omega$,
\begin{equation}
v \cdot n = 0,\label{L}
\end{equation}
meaning that the fluid particles cannot cross the boundary  ($n$ denotes the unit outward
normal). We suppose that the curl of \vinitial\ is a vortex patch, which involves some conormal
smoothness implying $\vinitial \in \Lip(\Omega)$ but not $\vinitial \in \cup_{\epsilon>0} \,
C^{1+\epsilon} (\Omega)$, and examine the classical problems of the existence of a solution,
either locally or globally in time, and of the persistence of the initial regularity.

\subsection{Brief history of the problem}

Vortex patches are initially a two-dimensional problem.
In 1963, V.~I.~{Yudovich} proved in~\cite{yud2} the existence and uniqueness of a weak, 2-D
solution of the Euler equations when the curl of \vinitial\ is a bounded function with compact
support; his result is valid even if the vorticity is not continuous. For example, the initial
vorticity may be the characteristic function of a bounded domain, or such a function multiplied
by a constant (that is, a so-called \emph{vortex patch}). The corresponding solution is not
Lipschitzian, but quasi-Lipschitzian in the sense that, $\forall t$,
\[
|v (t,x) - v (t,y)| \leq C |x-y|\ln (e+|x-y|)
\]
for a constant $C$ depending only on $t$. Also, $v$ has a flow $\psi$ which is, among other
properties, bicontinuous. In the case of a vortex patch, this implies that the vorticity
$\omega$ at time $t$ is the characteristic function of a domain which is homeomorphic to the
initial domain, because one has, in the two-dimensional case, $\omega(t, \psi(t,x)) = \omega_0
(x)$, $\forall t$, $\forall x\in\Omega$.

A more intricate question is how the smoothness of the \emph{boundary} of the patch evolves. In
1986, A.~{Majda} first conjectured, in view of numerical evidences, that the boundary of some
patches, initially regular, eventually produced singularities (see~\cite{majda}). But in 1991,
J.-Y.~{Chemin} proved the opposite result. The proof, detailed in~\cite{chemin}, makes use of
the notion of tangential smoothness along a system of vector fields: if a vortex patch has a
smooth boundary, say of class $C^{1+r}$, with $0<r<1$, one can construct vector fields of class
$C^r$ which are tangent to the boundary of the patch; the derivatives of the vorticity along
these vector fields, in the sense of distributions, have some regularity which is preserved by
the Euler system up to an arbitrary time.

Using similar techniques further results have been obtained, so today we are not restricted any
more to flows in the whole of $\rdeux$. Most important to us have been the generalization to
the dimension three by P.~{Gamblin} and X.~{Saint-Raymond} (see~\cite{gamstr}) and the article
of N.~{Depauw}, \cite{depauwarticle}, devoted to the problem of vortex patches in a bounded
domain of $\rdeux$. Numerous historical remarks and a comprehensive bibliography on the subject
can be found in~\cite{chemin} or in its English translation, \cite{chemine}.

\subsection{Notations}
We write as in~\cite{boubre} that
\[\domec=\{x\in\rtrois;\delta(x)>0\},\] and
\[\bord{\domec}=\{x\in\rtrois;\delta(x)=0\},\]
where $\delta$ is a \cinf\ function such that $n = -\nabla\delta$.

We use the notations of~\cite{chemin} for paradifferential calculus: $T$ is the paraproduct,
$R$ the remainder, $\tilde{h }$, $h \in \schw(\rtrois)$ have Fourier transforms $\chi$, $\varphi$
whose support are contained, respectively, in a ball and an annulus centered at the origin,
$\deln = \varphi(2^{-n}D)$ for $n \geq 0$, $\Delta_{-1} = \chi(D)$, etc.

Besov spaces on \rtrois\ are defined (see~\cite{triebel1}, section~2.3) by
\[
 \besovspq(\rtrois)
    = \{ f \in \schw'(\rtrois); \|f \|_{\besovspq(\rtrois)} \egnot \| \{2 ^{ns} \nLp{\deln
f}\}_{n=-1}^\infty \|_{l^q}  )
    < \infty \}
\]
and $\besovspq(\Omega)$ is the set of restrictions to $\Omega$ of all elements of
$\besovspq(\rtrois)$, so there is no distinction between $\besovspq(\Omega)$ and
$\besovspq(\ferm{\domec})$. Actually, we will only use $B^s_{1,2}$, in section~\ref{multipl},
and H{\"o}lder spaces $C ^r _*$ (or simply $C^r$ if $r \not\in{\mathbb Z}$), which correspond
to $p=q=\infty$.

Finally, we will note
\[
\Linf(X;\Cr(\domec))=\{w\in\bigcap_{r'<r}C(X;\Crp(\domec));\exists B\in\rr^+
:\nr{w(t)}\leq B,\ \forall t\in X\},
\]
for any interval (of time) $X$, and, $\forall r \in\,]0,1[$,
\[
[f]^\domec _r = \sup_{\substack{ x, y \,\in\, \domec\\ x \neq y}} \frac{|f (x) - f (y)| }{|x-y|
}.
\]

By now $r$ will denote a fixed real number in $]0,1[$.

\subsection{Results}

First of all we need to define what we mean by a ``vortex patch''.
\begin{definition}
We call vortex patch (or more precisely \Cr{} vortex patch) any vector field of the form
\[
\omegazero = \restr{ (\omega_{0i}\fcaract{\poche } +
\omega_{0e}\fcaract{\Omega\setminus\fermeture{\poche}}) }{\Omega},
\]
where $\omega_{0i}$, $\omega_{0e}\in C^r(\rtrois)$ and $\poche\subset\rtrois$ (the support of the patch) is an
open set of class \cunplusr.
\end{definition}
The following definition will also prove to be convenient.
\begin{definition}
Let $\poche \subset \rtrois$ be an open set of class \cunplusr. A system of vector fields $W=\{w ^{\nu};
\nu=1,\ldots,N'\}$ is \poche-regular if and only if
\begin{enumerate}
\item for some $s \in\,]0,r]$, $w ^{\nu}\in C ^s (\rtrois;\rtrois)$, $\forall\nu\in\{1,\ldots,
N'\}$,
\item each $w ^\nu$ is tangent both to \bord{\domec} and \bord{\poche}, and
\item $W$ is admissible, in the sense that (see~\cite[page 394]{gamstr})
\begin{equation}
[W]^{-1} \egdef \left\{\frac{2 }{N'(N'-1) }\sum_{\mu<\nu} |w^ \mu \produitvectoriel w^ \nu|
^2\right\}^{-1/4},
 \label{defW}
\end{equation}
is bounded on $\Omega$.
\end{enumerate}
\end{definition}
Now we can state our main result.
\begin{thm}
\label{thm1}
Let $v_0$ be a divergence-free vector field, tangent to $\bord{\Omega}$, whose curl is a \Cr{}
vortex patch, of support \poche. Suppose there exists a \poche-regular system of $C^s$ vector
fields, for some real $s\!\in \,]0,r]$. Then the Euler equations (\ref{E}), (\ref{I}),
(\ref{L}) have on $[0,T]$, for a time $T>0$, a (unique) solution $v \in \Linf ([0,T];
\Lip(\Omega))$.
Moreover, $\omega(t) \egnot \rotationnel v(t)$ remains a vortex patch, whose support
$\psi(t,P)$ is of class $C^{1+s}$, $\forall t \in [0,T]$, $\psi$ denoting the flow of $v$.
\end{thm}
Remark that, when $\fermeture{\poche} \subset \domec$, a \poche-regular system of $C ^r$ vector
fields can easily be constructed using proposition~3.2 in~\cite{gamstr} (page 395) and cut-off
functions. So our result is complete in that case. When \poche{} is tangent to $\bord{\domec}$,
however, theorem~\ref{thm1} doesn't always apply and, even if it does, there is a loss of
regularity. In this respect, theorem~\ref{thm1} is by far not as satisfactory as Depauw's
result in the 2-D case (namely, local existence of vortex patches tangent to the boundary and
preservation of full regularity, no matter how $\poche$ and $\bord{\domec}$ are tangent). But
our method also yields the following global results.
\begin{thm}
\label{thm2}
In addition of the hypotheses of theorem~\ref{thm1}, suppose that $v_0$ is two-dimensional or
axisymmetric, and in the later case, that $(\rotationnel v_0)/ \delta \in\Linf(\domec)$, where
$\delta$ is the distance to the axis of symmetry. Then the existence and regularity results are
in fact global in time, i.e.\ $T$ in the conclusion of theorem~\ref{thm1} can be taken
arbitrarily large.
\end{thm}
In the 2-D case, this shows that when a \poche-regular system of $C ^s$ vector fields exists,
Depauw's result of local existence for tangent patches can be completed by a global one (we
still can't rule out the possibility of a blow-up for $C^{1+r}$ norms, but the patch will
remain of class $C^{1+s}$ for all time). The difference of results between the two methods, in
the 2-D case, is mainly due to the estimate (\ref{dyna4}), which is better but less general
than the corresponding one in~\cite{depauwarticle} (indeed, it is the same estimate that one
can get when there is no boundary).

\subsection{Plan of the article}

We split the proof of theorem~\ref{thm1} in two sections. The existence of the solution is
shown first, in section~\ref{exists}. Essentially we adapt the methods of~\cite{gamstr}, using
extension procedures described in section~\ref{prolong}. The fact that $\omega(t)$ remains a
$C^s$ vortex patch is proved in section~\ref{regular}. Finally, in section~\ref{twodim}, we
explain how to get the global results in the 2-D and axisymmetric cases.

\section{Existence of the solution}
\label{exists}
\subsection{Sketch of the proof}\label{sketch}
The first part of theorem~\ref{thm1} is a consequence of the following result about
persistence of conormal smoothness.
\begin{prop} \label{propregtang}
Let $v _0$ be a divergence-free vector field, tangent to $\bord{\Omega}$, with $\omega_0
\egnot \rotationnel v_0 \in\Linf(\domec)$. Let also $W_0=\{w ^\nu_0; \nu=1,\ldots, N' \}$ be
an admissible system of $C^r$ vector fields that, again, are tangent to $\bord{\Omega}$.
Suppose that $\teste{\nabla}{w ^\nu_0 \otimes \omega_0 } \in C^{r-1}(\domec)$, $\forall \nu$,
and that $\omega_0 \cdot \tilde{n } \in C^r(\domec)$ for a field $\tilde{n}\in C^r(\rtrois)$
equal on $\bord{\Omega}$ to the unit outward normal to $\Omega$.

Then there exist a time $T_0>0$ and a constant $C_0$, both depending only on $
\normeLinf{\omega_0}$, $\| \omega_0 \cdot \tilde{n} \|_{C^r(\Omega)}$, $\normeLinf{[W_0] ^{-1
}}$, $ \normeCr{w ^\nu_0}$ and $\normeCrmu{\teste{\nabla}{w ^\nu_0 \otimes \omega_0 }}$,
$\nu=1,\ldots, N'$, such that the Euler equations (\ref{E}), (\ref{I}), (\ref{L}) have a
(unique) solution $v \in \Linf ([0,T_0]; \Lip(\Omega))$ and that $\normeLip{v(t,\cdot)} \leq
C_0$, $\forall t \in[0,T_0]$.
\end{prop}
Indeed, assume that proposition~\ref{propregtang} is valid. Then let $\omega_0$ be a vortex
patch of support \poche{} and $W_0$ be a \poche-regular system of $C ^s$ vector fields, as in
theorem~\ref{thm1}. The \poche-regularity of $W_0$ ensure the existence of a field $\tilde{n
} \in C ^s(\rtrois)$ equal on $\bord{P}$ to the unit outward normal to $P$ and equal on
$\bord{\Omega}$ to the unit outward normal to $\Omega$. As \omegazero{} is divergence-free,
$\omegazero\cdot \tilde{n }$ must be continuous (hence $C^s$) on $\Omega$.  Finally, the
proposition below expresses that \omegazero{} has good derivatives in the directions that are
tangent to the boundary of the patch.
\begin{prop} \label{typepoche}
With the above notations, let $w \in C^s(\rtrois)$ be a vector field which is tangent to
$\bord{P}$. The field $\teste{\nabla}{w \otimes ( \omega_{0i}\fcaract{P } + \omega_{0e}
\fcaract{\rtrois \setminus \fermeture{P}} ) }$ belongs to $C^{s-1}(\rtrois)$, and one has the
estimate
\begin{eqnarray*}
\lefteqn{ \normeCsmu{\teste{\nabla}{w \otimes ( \omega_{0i}\fcaract{P } + \omega_{0e}
\fcaract{\rtrois \setminus \fermeture{P}} )} } }\\
& \leq & C \normeLinf{w}(\normeCs{\omega_{0i}}+\normeCs{\omega_{0e}})+C
\normeLinf{\omegazero}\normeCs{w},
\end{eqnarray*}
for a constant $C$ depending only on the norm of the multiplication by $\fcaract{P}$ in
${\cal L} (C^{s-1}(\rtrois))$ (see section~\ref{multipl}).
\end{prop}
\demonstration
By regularization and application of the Gauss-Green formula, it is easily seen that
$\divergence(w\fcaract{P}) = \fcaract{P}\,\divergence w$ --- thanks to the fact that $w$ is
tangent to $\bord{P}$. Set $a = w\fcaract{P}$ and $b = \omega_{0i}$. We decompose again in
paraproducts and  remainders:
\begin{eqnarray*}
\teste{\nabla}{w \otimes \omega_{0i}\fcaract{P}} & = & \teste{\nabla}{a \otimes b} \\
& = & \partial_i[\paraproduit{a ^i}{b }+ \paraproduit{b }{a ^i}+ \reste{a ^i}{b }],
\end{eqnarray*}
and write the usual inequalities (see~\cite{chemin}, chapter~2),
\[
\normeCs{\paraproduit{a ^i}{b }+\reste{a ^i}{b }} \leq C \normeLinf{a ^i} \normeCs{b },
\]
and
\begin{eqnarray*}
\normeCsmu{\partial_i\paraproduit{b }{a ^i}} & \leq & \normeCsmu{\paraproduit{\partial_i b}{a
^i} } + \normeCsmu{\paraproduit{b }{\,\divergence a } } \\
& \leq & C \normeCsmu{\partial_i b } \normezer{a }+ C \normeLinf{b } \normeCsmu{\divergence a
}, \end{eqnarray*}
hence, by proposition~\ref{propmultipl} (in section~\ref{multipl},
page~\pageref{propmultipl}),
\[
\normeCsmu{\teste{\nabla}{w \otimes \omega_{0i}\fcaract{P}}} \leq C \normeLinf{w
}\normeCs{\omega_{0i }} + C \normeLinf{\omega_{0i }} \normeCs{w }.
\]
Now it suffices to do the same calculation with $b = \omega_{0e }$, using that
\[
\divergence(w \fcaract{\rtrois \setminus \fermeture{P}}) = \divergence w -
\fcaract{P}\,\divergence w.
\]
\findemonstration
So the first part of theorem~\ref{thm1} is a consequence of proposition~\ref{propregtang},
with $s$ substituted to $r$.

Now let us prove proposition~\ref{propregtang}. Assume first that $v _0
\in\cinf(\fermeture{\Omega})$. Then there exists a smooth solution $v$ of the Euler equations
defined on a time-interval $[0,T]$, with $T>0$; we just need to estimate $\normeLip{v(t)}$ on
any interval $[0,T]$ of existence of the  solution, independently of $T \leq T_0$, with $T_0$
to be determined.

Since the curl of $v$ is a solution of the system
\begin{equation}
\left\{\begin{array}{l}
\partial_t\omega+ \teste{\nabla}{v \otimes \omega} = \omega\cdot \nabla v, \\
\restr{\omega}{t=0 } = \omegazero,
\end{array}\right.                 \label{eqtourb}
\end{equation}
we have
\[
\normeLinf{\omega(t)} \leq \normeLinf{\omega_0(t)} e ^{\int _0^t\normeLip{v (s)}\,ds},
\]
and of course, if we define $w ^\nu$, $\nu=1,\ldots,N'$, as the solutions of
\begin{equation}
\left\{\begin{array}{l}
        \partial_t w ^\nu+ \teste{\nabla}{v \otimes w ^\nu} = w ^\nu \cdot \nabla v, \\
        \restr{w ^\nu}{t=0 } = w ^\nu_0,  \label{defwnu}
\end{array}\right.
\end{equation}
a similar inequality is true for each $\normeLinf{w ^\nu(t)}$. Set $W(t) = \{w ^\nu(t);
\nu=1,\ldots,N'\}$. The proof of corollary~4.3 in~\cite{gamstr} shows that
\[
\normeLinf{[W (t)]^{-1} } \leq C \normeLinf{[W _0] ^{-1 }} e ^{\frac{1 }{2 }
\int_0^t\normeLip{v (s)}\,ds},
\]
hence $W(t)$ is admissible. Moreover, \label{wnutangents}the fields $w ^\nu (t)$ are,
$\forall t \in[0,T]$, tangent to $\bord{\Omega}$. Indeed, it follows from the definition
of the $w ^\nu$'s, (\ref{defwnu}), that $w \pgrad (\delta \rond \psi^{-1 })$ is constant
along the particle trajectories
(see~\cite{majda}, proposition~2.1, page~S201; $\psi$ is the flow of $v$), thus $\restr{w
^\nu (t)}{\bord{\Omega}} \cdot n = 0$, $\forall t$. So we may apply, $\forall t$,
lemma~\ref{lemstatique} (see section~\ref{chapstatique}), which gives
\begin{equation}
\normeLip{v (t)} \leq C (1+\normeLinf{\omega(t)}) \ln (e+ X(t)), \label{eststatinstantt}
\end{equation}
with
\begin{multline*}
 X(t) \egnot 1  + \| \omega(t) \|_{\Linf(\Omega)}+ \| [W(t)]^{-1 } \|_{\Linf(\Omega)}
               \\
            +  \sum_{\nu=1 }^{N'} \| w ^{\nu}(t) \|_{C^r(\Omega)}
           + \sum_{\nu=1 }^{N'} \| \teste{\nabla }{ w ^{\nu}(t) \otimes \omega(t) }
\|_{C^{r-1}(\Omega)}
              + \| \omega(t)\cdot n \|_{C ^{r}(\bord{\Omega})}.
\end{multline*}
Remark that the terms in $X(t)$ are analogous to those found in~\cite{gamstr}: only
$\|\omega(t)\cdot n\|_{C ^r(\bord{\domec})}$ is a new one.
Using (\ref{defwnu}) and (\ref{eqtourb}), we show in section~\ref{dynam} the estimates
\begin{equation}\begin{split}                    \label{dyna13}
\normeCr{w ^\nu(t)}
& \leq C_0 e ^{C\int_0^t \normeLip{v (s)}\,ds}
+  C \int _0^t  \normeLip{v (s)} (\normeCr{w^\nu(s)}
\\  & + \frac{\normeCrmu{\teste{\nabla}{w ^\nu(s) \otimes \omega(s)}}
}{\normeLinf{\omega(s)}}) \, e
^{C\int_s^t \normeLip{v (\tau)}\,d \tau} ds,
\end{split}\end{equation}
\begin{multline}                                 \label{dyna15}
\frac{\normeCrmu{\teste{\nabla}{w ^{\nu}(t) \otimes \omega(t)}} }{\normeLinf{\omega(t)}}
\leq  C_0 e ^{C\int_0^t \normeLip{v (s)}\,ds} \\
+C \int_0^t \normeLip{v (s)} ( \normeCr{w ^\nu (s)}+\frac{\normeCrmu{\teste{\nabla}{w
^{\nu}(s) \otimes \omega(s)}} }{
\normeLinf{\omega(s)}} \\
+ \frac{\normeLinf{w ^\nu (s)}}{\normeLinf{\omega(s)}} \| \omega(s)\cdot n
\|_{C^r(\bord{\Omega})}) \, e ^{C \int_s^t \normeLip{v (\tau)}\, d \tau } ds,
\end{multline}
and
\begin{equation}                                  \label{dyna18}
\| \omega(t)\cdot n \|_{C^r(\bord{\Omega})} \\
 \leq  C _0 (1+ \sum_{\nu=1 }^{N' }  \normeCr{w ^\nu(t)}) \, e ^{C \int_0^t \normeLip{v
(s)}\,ds },
\end{equation}
$C_0$ denoting, from now on, any constant which depends only on the quantities
$\normeLinf{\omega_0}$, $\|\omega_0 \cdot n\|_{C^r(\bord{\Omega})}$, $\normeLinf{[W_0] ^{-1
}}$, $ \normeCr{w ^\nu_0}$ and $\normeCrmu{\teste{\nabla}{w ^\nu_0 \otimes \omega_0 }}$,
$\nu=1,\ldots, N'$. Summing up, we get
\begin{align*}
X (t) & \leq  C_0 e ^{C\int_0^t \normeLip{v (s)} } \,ds \\
&  + C_0 \int_0^t \normeLip{v (s)} e ^{C \int_0^s \normeLip{v (\tau)} \,d\tau} X (s) \, e ^{C
\int_s^t \normeLip{v (\tau)} \,d\tau} \,ds.
\end{align*}
Dividing by $e ^{C\int_0^t \normeLip{v (s)} } \,ds$ and applying Gronwall's lemma gives
\begin{equation}
X (t) \leq C_0 e ^{C_0 \int_0^t \normeLip{v (s)} e ^{C\int_0 ^s \normeLip{v (\tau)} \,d \tau}
\,ds,} \quad \forall t \in [0,T];\label{majxot}
\end{equation}
introducing this in (\ref{eststatinstantt}), we have
\begin{equation} \label{Deltath}
\begin{split}
\normeLip{v (t)} & \leq C (1 + \normeLinf{\omegazero} e ^{C\int_0^t \normeLip{v (s)}\,ds}) \\
& \quad \ln( e + C_0 e ^{C_0 \int_0^t \normeLip{v (s)} e ^{C\int_0^s \normeLip{v(\tau)} \, d
\tau } \,ds}),
\end{split}
\end{equation}
so, as $\ln(e+ a e ^b) \leq b + \ln(e+a)$, $\forall a,b \in \rr^+$,
\[
\normeLip{v (t)} \leq  C_0 e ^{C\int_0^t \normeLip{v (s)}\,ds},
\]
hence
\[
\normeLip{v (t)} \leq 2 C_0, \quad \forall t \in [0,\frac{1 }{2C C_0 }],
\]
following classical arguments: see for example~\cite[p.~407]{gamstr} (corollary~4.4). This
completes the proof of proposition~\ref{propregtang} if $v _0 \in\cinf(\fermeture{\Omega})$.

If $v _0 \not \in \cinf(\fermeture{\Omega})$, then we regularize $\widetilde{\omegazero}
\egnot \omega_{0i}\fcaract{P } + \omega_{0e} \fcaract{\rtrois \setminus \fermeture{P}}$ in
the usual way and we restrict the regularized fields to $\Omega$: \label{regularisation}
\[
\tilde{\omega}_{0k } \egdef \conv{\rho_{1/k}}{\widetilde{\omegazero}}, \quad \omega_{0k }
\egdef \restr{\tilde{\omega}_{0k}}{\Omega}, \quad \forall k \in {\mathbb N}_0.
\]
Let $v_k$, $ k \in \naturels_0$, be the smooth solutions of the Euler system corresponding to
the initial vorticities $\omega_{0k}$. It is clear that
$
\|\omega_{0k}\|_{\Linf(\Omega)}
%\leq \normeLinf{\tilde{\omega}_{0k}}
\leq C \normeLinf{\omegazero},
$
the constant $C$ being independent of $k$. We also have, $\forall \nu \in\{1,\ldots,N' \}$,
\begin{eqnarray*}
\lefteqn{ \|\teste{\nabla}{w ^\nu_0 \otimes \omega_{0k}}\|_{C^{r-1}(\Omega)} }\\
& \leq & \normeCrmu{\teste{\nabla}{w ^\nu_0 \otimes \tilde{\omega}_{0k}}} \\
& \leq & \normeCrmu{\teste{\nabla}{w ^\nu_0 \otimes \tilde{\omega}_{0k} -
\conv{\rho_{1/k}}{(w^\nu_0 \otimes \widetilde{\omegazero})}}} +
\normeCrmu{\conv{\rho_{1/k}}{\teste{\nabla}{w ^\nu_0 \otimes \widetilde{\omegazero}}}} \\
& \leq & C \normeCr{w ^\nu_0} \normeLinf{\widetilde{\omegazero}} + C
\normeCrmu{\teste{\nabla}{w^\nu_0 \otimes \widetilde{\omegazero}}} \\
& \leq & C \normeCr{w ^\nu_0} \normeLinf{\widetilde{\omegazero}} + C \normeLinf{w^\nu_0}
(\normeCr{\omega_{0i}}+\normeCr{\omega_{0e}}).
\end{eqnarray*}
and similarly
\begin{eqnarray*}
\|\omega_{0k} \cdot n\|_{C^s(\bord{\Omega})}
& \leq & \|\tilde{\omega}_{0k} \cdot \tilde{n}\|_{s} \\
%& \leq & \|\tilde{\omega}_{0k} \cdot \tilde{n} - \conv{\rho_{1/k}}{(\widetilde{\omegazero}
%\cdot \tilde{n})}\|_{s} + \|\conv{\rho_{1/k}}{(\widetilde{\omegazero} \cdot
%\tilde{n})}\|_{s}\\
%& \leq & C \|\tilde{n}\|_{s} \normeLinf{\tilde{\omega}_{0k}} + \|\widetilde{\omegazero}
%\cdot \tilde{n}\|_{s} \\
& \leq & C \|\tilde{n}\|_{s} \normeLinf{\omegazero} + C\|\widetilde{\omegazero} \cdot
\tilde{n}\|_{s}.
\end{eqnarray*}
Since proposition~\ref{propregtang} is already proved for regular initial data, these uniform
estimates imply that the $v _k$'s are defined on a common interval $[0,T_0]$ for some time
$T_0>0$, with the existence of a constant $C_0$ such that
\begin{equation}
\normeLip{v_k(t,\cdot)} \leq C_0, \quad \forall t \in [0,T_0],\quad \forall k \in
\naturels_0. \label{theo1}
\end{equation}
Using that the $v _k$'s are solutions of the Euler system, one can derive from (\ref{theo1})
a uniform bound on $\normeLinf{\partial_t v_k}$. So there exists
%one can therefore conclude by a compactness argument. One has, $\forall k \in\naturels$,
%
%\begin{eqnarray}
%\normeLinf{\partial_t v_k} & \leq & \normeLinf{v_k} \normeLip{v_k} + \|\nabla\pi(v_k,v
%_k)\|_{\sigma} \nonumber\\
%& \leq & \normeLinf{v_k} \normeLip{v_k} + C \|v_k\|_{\sigma} \normeLip{v_k} \nonumber\\
%& \leq & C, \label{theo2}
%\end{eqnarray}
%
%for a constant $C$ independent of $k$ and any $\sigma$ in $]0,1[$, by proposition~1.5
%in~\cite{mathese} (p.\ 19). (Here $\nabla\pi$ is a bilinear operator defined as follows. Let
%$v$ and $w$ be two vector fields of class $C^1$ on
%$\ferm{\domec}$. Then $\nabla\pi(v,w)$ is the gradient of any solution, in the Sobolev space
%$\sobud(\domec)$, of the Neumann problem
%
%\[
%\left\{\begin{array}{l}
%\lap\pi=-\sumijd\diwj\djvi\ \mbox{\ \ \ in \domec}\\
%\dsdn{\pi}=-\sumijd\wi\vj\di\djj\delta+C(v,w)\ \mbox{\ \ \ on\ }\bord{\domec}
%\end{array}\right.
%\]
%
%where
%
%\[C(v,w)\egnot\frac{1}{\meas_{d-1}(\bord{\domec})}\sumijd(-\int_\domec\diwj\djvi+
%\int_{\bord{\domec}}\wi\vj\di\djj\delta).)\]
%
%From (\ref{theo1}) and (\ref{theo2}), one deduces the existence of
a subsequence of $(v _k)$, $(v_{k(l)})_{l \in {\mathbb N}}$ converging in $\Linf([0,T_0];
C^\sigma(\Omega))$, $ \forall \sigma\in\,]0,1[$, to a field $v \in \Linf([0,T_0];
\Lip(\Omega))$. Then, since $v_k(0) \rightarrow v_0$ in $\Linf(\Omega)$,
%As
%\begin{eqnarray*}
%\|\nabla\pi(v _k,v _k) - \nabla\pi(v,v)\|_{\sigma} & = & \|\nabla\pi(v_k-v,v_k+v)\|_{\sigma}
%\\
%& \leq & C \|v _k- v \|_{\sigma} \normeLip{v _k + v},
%\end{eqnarray*}
%always by proposition~1.5 in~\cite{mathese},
it is easily seen that $v$ is the desired solution.

\subsection{Extensions of divergence-free vector fields}
\label{prolong}

Before we can prove the estimates (\ref{eststatinstantt}), (\ref{dyna13}), (\ref{dyna15}) and
(\ref{dyna18}), we must present some new extension lemma's. The method which is normally used
to extend divergence-free vector fields (see~\cite{giraultraviart}, for example) would not,
unfortunately, be sufficient for our purpose.

These lemma's rest upon the following elementary construction.

\subsubsection{Basic construction}\label{bc}

At every point $x\in\bord{\Omega}$ (here $\Omega$ could be any bounded, \cun\ open subset of
$\rtrois$), we choose orthogonal vectors, $e^x_1$, $e^x_2$ and $e^x_3$, each of length one, which are
not tangent to $\bord{\Omega}$ and are oriented towards the exterior of $\Omega$.

It is not hard to check that, $\forall x\in \bord{\Omega}$, there exist two neighbourhoods of $x$,
$W_x$ and $V_x \fermeturedans W_x$, and Lipschitzian projections $y^{x,j}: V_x \mapsto W_x$ such that
$\forall \xi \in \fermeture{V_x}$, $ \forall j \in \{1,\,2,\,3\}$, the line through $\xi$ and of
direction $e^x_j$ meets $\bord{\Omega}\cap W_x$ in exactly one point, $y^{x,j}(\xi)$.

The collection of all $V_x$ is an open cover of the compact set $\bord{\Omega}$. So we can extract from
it a finite subcover $V_1,\ldots,V_N$ (each $V_i$ corresponding to a point $x_i \in \bord{\Omega}$) and
after that choose open sets $V_0$ and $V$ such that $\fermeture{V_0} \subset \Omega \subset
\bigcup_{i=0}^N V_i \fermeturedans V$. Finally we give ourselves $N+1$ functions $\psi_i \in
\cinfc(\rtrois)$ ($i=0,\ldots,N$), with values in $[0,1]$, such that $\support{\psi_i} \fermeturedans
V_i$ and $\sum_{i=0}^N \psi_i = 1$ on $ \fermeture{\Omega}$.

\subsubsection{Extensions}
In this section we denote by $\uronde$ the set of continuous divergence-free vector fields on
$\fermeture{\domec}$, and by $\cronde$ the set of all continuous vector fields on \rtrois.

Our first extension operator, $P$, maps elements of \uronde{} to fields of bounded divergence.
\begin{lem} \label{lem1}
There exists an operator $P : \uronde \rightarrow \Linf(\rtrois; \rtrois)$ such that, $\forall u
\in\uronde$, $Pu$ extends $\restr{u}{\Omega}$, $\support{Pu} \subset V$ and
\begin{subequations}
\begin{equation}
\| Pu \|_{\Linf(\rtrois)} \leq C \| u \|_{\Linf(\Omega)}, \label{lem1a}
\end{equation}
\begin{equation}
\| Pu \|_{C^r(\rtrois \setminus\fermeture{\Omega})} \leq C \| u\cdot n \|_{C^r(\bord{\Omega})},
\label{lem1b}
\end{equation}
\begin{equation}
\| \divergence \tilde{u}\|_{\Linf(\rtrois)} \leq C \| u \|_{\Linf(\Omega)}, \label{lem1c}
\end{equation}
\begin{equation}
\| \divergence \tilde{u} \|_{C^r(\rtrois\setminus\fermeture{\Omega})} \leq C \| u\cdot n
\|_{C^r(\bord{\Omega})}, \label{lem1d}
\end{equation}
\end{subequations}
where $C$ depends only on \domec.
\end{lem}
\demonstration
On every $V_i$ ($i=1,\ldots,N$), we define a new vector field $\tilde{u}_i$ by
\begin{equation}
\label{utildeidef}
\tilde{u }_i(x) =
\begin{cases}
u (x) & \text{if }x \in V_i \cap\Omega, \\
[u(x) \cdot n(x)] \, n (x) & \text{if }x \in V_i \cap\bord{\Omega}, \\
\sum_{j=1}^{3 } [\tilde{u }_i(y^{i,j}(x)) \cdot e^{x_i}_j] \, e^{x_i}_j & \text{if }x \in V_i \setminus
\fermeture{\Omega},
\end{cases}
\end{equation}
with $y^{i,j} \egnot y^{x _i,j}$. Then we set
\begin{equation}
Pu = \psi_0 u + \sum_{i=1}^N \psi_i \tilde{u}_i, \label{prol3}
\end{equation}
so the estimate (\ref{lem1a}) is trivial. Let us prove (\ref{lem1b}).
We have $\| \tilde{u}_i \|_{\Linf(V _i\setminus\Omega)} \leq \| u\cdot n \|_{ \Linf(\bord{\Omega})}$
and if $\xi,\,\eta \in V _i\setminus\Omega$, then
\begin{eqnarray}
\lefteqn{ | \tilde{u}_i(\xi)-\tilde{u}_i(\eta) | } \nonumber \\
& \leq & \sum_{j=1}^3 | \tilde{u}_i(y^{i,j}(\xi))-\tilde{u}_i(y^{i,j}(\eta)) |\nonumber\\
& = & |[(u\cdot n)(y^{i,j}(\xi))] n(y^{i,j}(\xi)) - [(u\cdot n)(y^{i,j}(\eta))] n(y^{i,j}(\eta))|
:\nonumber
\end{eqnarray}
as $n$ is smooth, we only have to remember that the mappings $y^{i,j}$ are Lipschitzian.
In order to prove (\ref{lem1c}), remark that, $\forall i \in \{1,\ldots,N\}$,
$\divergence\tilde{u}_i=0$ on $V_i$, by construction. Thus
\[
\divergence \tilde{u} = u \cdot \nabla\psi_0+ \sum_{i=1}^N \tilde{u}_i\cdot \nabla\psi_i,
\]
and (\ref{lem1c}) holds with $C=  \sum_{i=0}^N \| \nabla\psi_i\|_{\Linf}$. The proof of (\ref{lem1d})
is analogous.
\findemonstration

Typically, we will apply lemma~\ref{lem1} to regular curls $\omega (t)$ of solutions corresponding to
regularized data. Then, introducing a discontinuity on \bord{\domec} (which may seem unnatural at
first) will be a crucial trick, because $\restr{\omega (t)\cdot n}{\bord{\domec}}$ is easier to control
than $\restr{\omega (t)}{\bord{\domec}}$ (see section~\ref{evolomn}), and also vanish in the 2-D and
axisymmetric cases.

Of course, we can also extend $u$ continuously through $\bord{\Omega}$, but then the tangential
component of $u$ on $\bord{\Omega}$ appears in the estimates. Doing this will be helpful only if $u
\in\Lip$.

\begin{lem} \label{lem2}
There exists an operator $P_c : \uronde \rightarrow \cronde$ such that, $\forall u \in\uronde$, $P_c u$
extends continuously $\restr{u}{\Omega}$, $\support{P_c u} \subset V$, and
%\begin{enumerate}
%\item if $u \in C^r(\Omega)$, then $\| P_c u \|_{C^r(\rtrois)} +\| \divergence P_c u
%\|_{C^r(\rtrois)} \leq C \| u \|_{C^r(\Omega)}$,
%\item
%if $u \in \Lip(\Omega)$, then $
\[
\| P_c u \|_{\Lip(\rtrois)} + \| \divergence P_c u \|_{\Lip(\rtrois)} \leq C \| u \|_{\Lip(\Omega)},
\]
%\end{enumerate}
where $C$ depends only on \domec.
\end{lem}
\demonstration Set, $\forall x \in V_i \cap \bord{\Omega}$, $\tilde{u}_i (x)= u (x)$ instead of $[u(x)
\cdot n(x)] \, n (x)$ in the definition of $\tilde{u}_i$, (\ref{utildeidef}), and define $P_c$ by the
right member of (\ref{prol3}).
\findemonstration

The third lemma can be used if one really wants the extended fields to be of free divergence. We note
$\urondep = \{ u \in\uronde; \int_B u \cdot n = 0$ for every connected component $B$ of
$\bord{\domec}\}$.
\begin{lem}\label{lem3}
There exists an operator $P_{\divergence}:\urondep\rightarrow\Linf(\rtrois;\rtrois)$ such that,
$\forall u \in\urondep$, $P_{\divergence} u$ extends $\restr{u }{\domec}$, $\support{P_{\divergence}u }
\subset\fermeture{V }$, $\divergence P_{\divergence} u =0$ on \rtrois, and
\begin{subequations}
\begin{equation}
\| P_{\divergence}u \|_{\Linf(\rtrois)} \leq C \| u \|_{\Linf(\Omega) }, \label{lem3a}
\end{equation}
\begin{equation}
\| P_{\divergence}u \|_{C^r(V \setminus \fermeture{\Omega})} \leq C \| u\cdot n
\|_{C^r(\bord{\Omega})}, \label{lem3b}
\end{equation}
\end{subequations}
the constant $C$ depending only of \domec.
\end{lem}
\demonstration Let us set $P_{\divergence} u = Pu -  \fonctioncaracteristique{V
\setminus\fermeture{\domec}}\nabla\psi$, where $\psi$ is a solution of
\[\left\{\begin{array}{l}
\lap\psi= \divergence P u \ \mbox{ in }V\setminus\ferm{\Omega} \\
\dsdn{\psi} = 0 \quad\mbox{ on }\bord{V}\cup\bord{\Omega}.
         \end{array}
  \right.\]
As $\| \nabla\psi \|_{C^1_\etoile (V \setminus \fermeture{\Omega})} \leq C \| \divergence  Pu
\|_{\Linf( V \setminus \fermeture{\Omega})}$, (\ref{lem3a}) and (\ref{lem3b}) are easy consequences of
(\ref{lem1a}) and (\ref{lem1b}). The other properties of $P_{\divergence}$ are obvious.
\findemonstration

\subsection{Static estimates}
\label{chapstatique}
Here we prove (\ref{eststatinstantt}), together with the estimate (\ref{estjbord}), which will be needed in
section~\ref{regular}.
\begin{lem} \label{lemstatique}
Suppose that $\omega= \rotationnel v$, where $v$ is a smooth (\cinf) divergence-free vector field such that $v
\cdot n =0$ on $\bord{\Omega}$. Let $W = \{w^\nu; \nu=1,\ldots,N'\}$, $N'\geq 2$ be an admissible family of
$C^r$ vector fields which, as $v$, are tangent to the boundary of $\Omega$. Then we have the estimate
\begin{equation}
\normeLip{v } \leq C (1+\normeLinf{\omega}) \ln (e+ X), \label{eststatique}
\end{equation}
with
\begin{multline*}
X \egnot 1  + \| \omega \|_{\Linf(\Omega)}+ \| [W]^{-1 } \|_{\Linf(\Omega)}
\\
            +  \sum_{\nu=1 }^{N'} \| w ^{\nu} \|_{C^r(\Omega)}
           + \sum_{\nu=1 }^{N'} \| \teste{\nabla }{ w ^{\nu} \otimes \omega } \|_{C^{r-1}(\Omega)}
              + \| \omega\cdot n \|_{C ^{r}(\bord{\Omega})},
%\label{defX}
\end{multline*}
and there exists a constant $C$ such that, for any subset $\Omega'$ of $\Omega$,
\begin{equation}
[\nabla v]^{\Omega'}_r \leq C X^{20} + [\omega]^{\Omega'}_r. \label{estjbord}
\end{equation}
\end{lem}
The proof of lemma~\ref{lemstatique} is spread over several subsections.
\subsubsection{Extensions}
A priori, the fields $\omega$ and $w ^\nu$, $\nu=1,\ldots,N'$, are only defined on $\Omega$. Careful
extensions, however, bring most problems back to the case of fields defined on the whole of \rtrois. We shall
use the notations of section~\ref{prolong}.

We simply extend $\omega$  to $ \bar{\omega} \egnot P_{\divergence}\omega$ (see lemma~\ref{lem3}).  The
extensions of the $w^ \nu$'s is more complicated.

Actually, every vector field $w^ \nu$ is extended $N$ times (see section~\ref{bc}). We define
$\tilde{w}^{\nu,i}$ ($\nu=1,\ldots,N'$, $i=1,\ldots,N$) on \rtrois\ in two steps:  we first set $\tilde{w
}^{\nu,i }= w ^{\nu}$ on \fermeture{\Omega} and $\tilde{w }^{\nu,i }= w ^{\nu}\rond y^{i,1 }$ (for example) on
\fermeture{V_i}; next we extend this auxiliary field on \rtrois, {\`{a}} la Whitney (see~\cite{stein},
chapter~VI, section~2). As these $\tilde{w }^{\nu,i }$'s have no reason to be tangent to \bord{V }, one cuts
them by a function  $\varphi_{\mbox{int}}\in\cinfc(\rtrois)$ equal to 1 on $\cup_{i=0 }^N V_i$ and whose
support is a subset of $V$. We note $\bar{w} ^{N(\nu-1)+i } = \varphi_{\mbox{int}} \tilde{w }^{\nu,i }$,
$\forall\nu,i$.

Of course, the system of the $\bar{w} ^{\mu}$'s, $\mu=1,\ldots,NN'$, is not admissible away from $\cup_{i=0
}^N V_i$; therefore we must add to this system other vector fields. By proposition~3.2 in~\cite{gamstr}, page
395, there exists an admissible system of vector fields $\{\tilde{w }^{NN'+j}; j=1,\ldots,5\}$, tangent to
\bord{V } and of class \cinf. These fields will not be tangent to \bord{\Omega}; so we set this time $\bar{w}
^{NN'+j } = \varphi _{\mbox{ext}}\tilde{w }^{NN'+j }$, $\forall j$, where $\varphi _{\mbox{ext}}
\in\cinf(\rtrois)$ is equal to 1 outside $\cup_{i=0 }^N V_i$ and vanish on a neighbourhood of
\fermeture{\Omega}. The $\bar{w}^{NN'+j }$'s depend only of $V$.

In that way we obtain vector fields $\bar{w }^{\mu}$ tangent to \bord{\Omega}\ and to \bord{V }, such that
\begin{equation}
\sum_{\mu=1 }^{NN'+5 } \| \bar{w }^ \mu \|_{C^r(\rtrois)} \leq C (1+ \sum_{\nu=1 }^{N' } \| w^ \nu
\|_{C^r(\Omega)}), \label{abc}
\end{equation}
and which form an admissible system:
\begin{equation}
[\bar{W }^{-1 }] \leq C(1+[W]^{-1 }), \label{def}
\end{equation}
because $[\bar{W }]^{-1 } \leq [W]^{-1 }$ on $\cup_{i=0 }^N V_i$ and $[\bar{W }]^{-1 } \leq C<\infty$ outside.

\subsubsection{Estimation of $\left\protect\langle\nabla, \,
\protect\bar{w}^{\mu}\otimes\protect\bar{\omega}\right\protect\rangle$ in $C^{r-1}(\rtrois)$}
\label{nablawbomb}

Let $\varphi\in\schw(\rtrois)$ be a test-function, and let $\bar{w }^ \mu_q = \chi(2^{-q}D) \bar{w }^ \mu$, $q
\in {\mathbb N}$, be a sequence of smooth fields obtained from $\bar{w}^ \mu$ ($\mu\in\{1,\ldots,NN'+5\}$ is
fixed) by regularization. We have, with summation on $i$,
\begin{eqnarray}
\lefteqn{ \teste{
\teste{\nabla}{\bar{w }^ \mu \otimes \bar{\omega}}
}{\varphi} } \nonumber \\
& = & - \lim_{q \rightarrow\infty } \int_{\Omega} \bar{w}^{\mu i}_q \bar{\omega} \,\partial_i \varphi
  - \lim_{q \rightarrow\infty } \int_{V \setminus \fermeture{\Omega}} \bar{w}^{\mu i}_q \bar{\omega}
\,\partial_i \varphi \nonumber\\
& = &  \lim_{q \rightarrow\infty } \int_{\Omega} \partial_i(\bar{w}^{\mu i}_q \bar{\omega}) \, \varphi
- \lim_{q \rightarrow\infty}\int_{\bord{\Omega}}\bar{w}^{\mu}_q\cdot n \,\trace\restr{\bar{\omega}}{\Omega}
\varphi
     \nonumber\\
& & \mbox{} + \lim_{q \rightarrow\infty } \int_{V \setminus \fermeture{\Omega}} \partial_i(\bar{w}^{\mu i}_q
\bar{\omega})  \varphi- \lim_{q \rightarrow\infty}\int_{\bord{\Omega}}\bar{w}^{\mu}_q\cdot (-n)
\,\trace\restr{\bar{\omega}}{V \setminus\fermeture{\Omega}} \varphi
         \nonumber\\
& & \mbox{}
   - \lim_{q \rightarrow\infty}\int_{\bord{V}}\bar{w}^{\mu}_q\cdot n \,\bar{\omega} \,\varphi, \label{stat1}
\end{eqnarray}
$n$ denoting again the unit outward normal, either to $\Omega$ or to $V$. When $q \rightarrow\infty$, the
integrals on \bord{\Omega}\ and on \bord{V}\ tend to 0, because $\bar{w}^{\mu}_q\cdot n \rightarrow
\bar{w}^{\mu}\cdot n$ $(= 0)$ uniformly, both on \bord{\Omega}\ and on \bord{V}. Now let $\tilde{P}$ be any
continuous extension operator, $\forall s \in\,]-2,2[$, from $C^s(\Omega)$ into $C^s(\rtrois)$
(see~\cite[section~3.3.4]{triebel1} for the existence of such operators). The distribution $\varphi \mapsto
\int_{\Omega} \partial_i(\bar{w}^{\mu i}_q \bar{\omega})  \varphi$ is equal to the usual product
$\fcaract{\Omega}\,\tilde{P}[\restr{\partial_i(\bar{w}^{\mu i}_q \bar{\omega})}{\Omega}]$.  As
$\tilde{P}[\restr{\partial_i(\bar{w}^{\mu i}_q \bar{\omega})}{\Omega}]$ $\rightarrow$
$\tilde{P}[\restr{\partial_i(\bar{w}^{\mu i} \bar{\omega})}{\Omega}]$, when $q \rightarrow\infty$, in $C
^{r'-1 }(\rtrois)$, $\forall r'<r$, the products $\fcaract{\Omega}\,\tilde{P}[\restr{\partial_i(\bar{w}^{\mu
i}_q \bar{\omega})}{\Omega}]$ converge to $\fcaract{\Omega}\,\tilde{P}[\restr{\partial_i(\bar{w}^{\mu i}
\bar{\omega})}{\Omega}]$ in the same spaces, by virtue of proposition~\ref{propmultipl}. So we have
\begin{equation}
\lim_{q \rightarrow\infty} \int_{\Omega} \partial_i(\bar{w}^{\mu i}_q \bar{\omega})  \varphi
= \teste{\fcaract{\Omega}\,\tilde{P}[\restr{\partial_i(\bar{w}^{\mu i} \bar{\omega})}{\Omega}]}{\varphi}.
\label{stat2}
\end{equation}
In the same way, one shows that
\begin{equation}
\lim_{q \rightarrow\infty } \int_{V \setminus \fermeture{\Omega}} \,\partial_i(\bar{w}^{\mu i}_q \bar{\omega})
\varphi
= \teste{\fcaract{{V \setminus \fermeture{\Omega}}}\,\partial_i[\tilde{P}\restr{(\bar{w}^{\mu i}
\bar{\omega})}{V \setminus\fermeture{\Omega}}\,]}{\varphi}.
\label{stat3}
\end{equation}
Putting (\ref{stat1}), (\ref{stat2}) and (\ref{stat3}) together, we get
\[
\teste{\nabla}{\bar{w}^{\mu}\otimes\bar{\omega}} = \fcaract{\Omega}\,\tilde{P}[\restr{\partial_i(\bar{w}^{\mu
i} \bar{\omega})}{\Omega}] + \fcaract{{V \setminus
\fermeture{\Omega}}}\,\partial_i[\tilde{P}\restr{(\bar{w}^{\mu i} \bar{\omega})}{V
\setminus\fermeture{\Omega}}\,],
\]
hence the estimate we wanted:
\begin{eqnarray*}
\lefteqn{ \sum_{\mu=1}^{NN'+5} \| \teste{\nabla}{\bar{w }^{\mu}\otimes\bar{\omega}} \|_{ C^{r-1}(\rtrois) }}
\nonumber \\
& \leq & C \sum_{\mu=1}^{NN'+5} \| \restr{\partial_i(\bar{w}^{\mu i}\bar{\omega})}{\Omega} \|_{C^{r-1
}(\Omega)} \\
& & \mbox{}+  \sum_{\mu=1}^{NN'+5} \sum_{i=1}^3 [\,\| \bar{w}^{\mu i} \|_{C^r(V \setminus \fermeture{\Omega})}
\| \bar{\omega} \|_{\Linf(V \setminus \fermeture{\Omega})} +  \| \bar{w}^{\mu i} \|_{\Linf(V \setminus
\fermeture{\Omega})} \| \bar{\omega} \|_{C^r(V \setminus \fermeture{\Omega})} ] \\
& \leq & C'\,[ \,\sum_{\nu=1 }^{N'}\| \teste{\nabla}{w^{\nu}\otimes\omega} \|_{C^{r-1 }(\Omega)} \\
& & \mbox{}+  (1+\sum_{\nu=1}^{N'} \| w^{\nu} \|_{C^r(\Omega)}) \| \omega \|_{\Linf(\Omega)} +
(1+\sum_{\nu=1}^{N'} \| w^{\nu} \|_{\Linf(\Omega)}) \| \omega\cdot n \|_{C^r(\bord{\Omega})}],
\end{eqnarray*}
thanks to (\ref{abc}), (\ref{lem3a}), (\ref{lem3b}) and proposition~\ref{propmultipl}.
\subsubsection{Biot-Savart's law}
To $\bar{\omega}$ we associate the field $\bar{v} \egdef \bar{\omega}\wedge \nabla F$, whose components are
\[
\left\{\begin{array}{l}
\bar{v}^1 = \conv{\bar{\omega}^2}{\partial_3 F }- \conv{\bar{\omega}^3}{\partial_2 F }, \\
\bar{v}^2 = \conv{\bar{\omega}^3}{\partial_1 F }- \conv{\bar{\omega}^1}{\partial_3 F }, \\
\bar{v}^3 = \conv{\bar{\omega}^1}{\partial_2 F }- \conv{\bar{\omega}^2}{\partial_1 F },
\end{array} \right.
\]
where $F(x)=-1/4 \pi |x|$ is the 3-D Laplacian's fundamental solution. We have identically $\divergence
\bar{v} = 0$, and $\rotationnel \bar{v} = \conv{\bar{\omega}}{\Delta F }
-\conv{\divergence\bar{\omega}}{\nabla F }  = \bar{\omega},$
but in general $\restr{\bar{v}\cdot n}{\bord{\Omega}} \neq 0$. So $\bar{v}$ is not an extension of $v$; as a
matter of fact, since $\Omega$ is simply connected, \label{pagequatre} $v= \restr{\bar{v}}{\Omega}-
\nabla\alpha$, where $\alpha$ is a solution of
\[
\left\{\begin{array}{l}
\Delta\alpha=0 \quad \mbox{in } \Omega,\\
\frac{\partial\alpha}{\partial n } = \bar{v}\cdot n \quad \mbox{on } \bord{\Omega}.
\end{array}\right.
\]

\subsubsection{Estimation of $v$ in $\Lip(\Omega)$}

Let us first estimate $\normeLip{\nabla\alpha}$. We use the inequality
\begin{equation}
\| \nabla\alpha \|_{\Lip} \leq C \| \nabla\alpha \|_{\Cunetoile(\Omega)} \ln (e +
                                                           \frac{ \| \nabla\alpha \|_{1+r}}
                                                                { \| \nabla\alpha \|_{\Cunetoile(\Omega)}}).
\label{stat4}
\end{equation}
We have directly
%\begin{eqnarray}
\[\| \nabla\alpha \|_{\Cunetoile(\Omega)}
%&
 \leq
 %&
 C \| \bar{v}\cdot n \|_{\Cunetoile(\bord{\Omega})} %\nonumber\\
%&
 \leq
 %&
 C \| \omega \|_{\Linf(\Omega)},\] %\label{stat5}
%\end{eqnarray}
and
%\begin{eqnarray}
\[\| \nabla\alpha \|_{1+r}
%&
\leq
%&
C \| \bar{v}\cdot n \|_{\cunplusr(\bord{\Omega})} %\nonumber\\
%&
\leq
%&
C \| \bar{v} \|_{\cunplusr(V \setminus \fermeture{\Omega})},\] %\label{stat6}
%\end{eqnarray}
by theorem~3.3.3 in~\cite{triebel1} (trace theorem).

To get the needed estimates on $\bar{v}$, we follow the proof of proposition~3.3 in~\cite{gamstr}. One makes
use of the Fourier multiplier $\Lambda= \lambda(D)$, where $\lambda(\xi) = [\chi(\xi)+|\xi|^2]^{1/2}$ (the
function $\chi\in\cinfc(\rtrois)$ is positive and equal to 1 near 0). As $\rotationnel\bar{v}=\bar{\omega}$,
\begin{eqnarray*}
\bar{v} & = & \Lambda^{-2}(\chi(D)-\Delta)\bar{v} \\
& = & \Lambda^{-2}\chi(D)\bar{v} + \Lambda^{-2}\rotationnel \bar{\omega} ;
\end{eqnarray*}
thus, $\forall i,j \in \{1,\,2,\,3\}$,
\begin{eqnarray}
\partial_j\bar{v}^i & = & \partial_j\Lambda^{-2}\chi(D)\bar{v}^i + \Lambda^{-2}
\partial_j(\rotationnel\bar{\omega})^i \nonumber\\
& = & \partial_j\Lambda^{-2}\chi(D)\bar{v}^i + \sum_{k,l=1}^3 \alpha_{k,l }^i \Lambda^{-2}
\partial_j\partial_k\bar{\omega}^l, \label{stat7}
\end{eqnarray}
with $\alpha_{k,l }^i\in\{-1,\,0,\,1\}$, $\forall i,k,l$. The functions $\partial_j\Lambda^{-2}\chi(D)\bar{v}$
($j=1,\,2,\,3$) are highly regular: $\| \chi(D) \bar{v} \|_{s } \leq C _s \normeLinf{\bar{v}} \leq C' _s \|
\omega \|_{\Linf(\Omega)}$, $\forall s \in\rr$. The difficult terms $\partial_j\partial_k\,\bar{\omega}^l$ are
treated as follows. By lemma's~3.4 and~3.5 in~\cite{gamstr}, there exists $C ^r$ functions $a _{jk }$ and $b
_{jk }^{l \nu}$ such that
\[
\xi_j\xi_k-a _{jk}|\xi|^2 = \sum_{l,\nu}b ^{l \nu}_{jk }\xi_l\teste{\bar{w }^{\nu}}{\xi},
\quad\forall(x,\,\xi)\in\rtrois\produitcartesien\rtrois,
\]
with the estimates
\[
\| a _{jk } \|_{r }\leq C (1+ \sum_{\mu=1 }^{NN'+5 }\normeLinf{\bar{w }^\mu})^7 (\sum_{\mu=1 }^{NN'+5 } \|
\bar{w }^\mu \|_{r })(1+ \normeLinf{[\bar{W}]^{-1 }}) ^{4 }
\]
and
\[
\| b _{jk }^{l \nu} \|_{r }\leq  C (1+ \sum_{\mu=1 }^{NN'+5 }\normeLinf{\bar{w }^\mu})^{18} (\sum_{\mu=1
}^{NN'+5 } \| \bar{w }^\mu \|_{r }) (1+ \normeLinf{[\bar{W}]^{-1 }}) ^{20 }.
\]
The $\xi$-identity implies
\[
\partial_j\partial_k\,\bar{\omega}- \Delta(a _{jk }\,\bar{\omega}) = \sum_{l,m,\nu}\partial_l\partial_m (b ^{l
\nu}_{jk }\bar{w} ^\nu_m \bar{\omega}),
\]
hence (see~corollary~3.6 in~\cite{gamstr}, pages 399 and 400)
\begin{eqnarray}
\lefteqn{ \| \partial_j\partial_k\,\bar{\omega}- \Delta(a _{jk }\,\bar{\omega}) \|_{r-2 }} \nonumber\\
& \leq & C \sum_{l,m=1}^{3 } \sum_{\nu=1}^{NN'+5 } ( \normeLinf{\bar{w} ^\nu_m} \normeLinf{\bar{\omega}} \| b
^{l \nu}_{jk } \|_{r } + \| \partial_m (\bar{w} ^\nu_m \bar{\omega}) \|_{r-1 } \normeLinf{b ^{l \nu}_{jk }})
\nonumber\\
& \leq & C X^{20}.
\end{eqnarray}
Writing
\begin{eqnarray}
\lefteqn{ \Lambda^{-2 }\partial_j\partial_k \,\bar{\omega}^l } \nonumber \\
& = & \Lambda ^{-2 }[\partial_j\partial_k \,\bar{\omega}^l - \Delta(a _{jk }\,\bar{\omega}^l) ] + \Lambda^{-2
}\Delta(a _{jk }\,\bar{\omega}^l) \nonumber\\
& = & \Lambda ^{-2 }[\partial_j\partial_k \,\bar{\omega}^l - \Delta(a _{jk }\,\bar{\omega}^l) ] +
\Lambda^{-2}\chi(D)(a _{jk }\,\bar{\omega}^l) - a _{jk }\,\bar{\omega}^l, \label{stat8}
\end{eqnarray}
it can be seen, firstly that
\[
\| \nabla\alpha \|_{1+r} \leq \| \nabla\otimes\bar{v } \|_{C^r(V \setminus \fermeture{\Omega})} \leq CX^{20},
\]
thus
\[
\| \partial_j\nabla\alpha \|_{\Linf} \leq C  \normeLinf{\omega}\ln(e+ \frac{X^{20} }{\normeLinf{\omega}}),
\]
and secondly that
\begin{eqnarray*}
\normeLinf{\partial_j\bar{v }} & \leq & C \normeLinf{\bar{\omega}} + \sum_{k,l =1}^3 \normeLinf{\Lambda
^{-2}[\partial_j\partial_k \,\bar{\omega}^l - \Delta(a _{jk }\,\bar{\omega}^l) ] } \\
& \leq & C \normeLinf{\omega} + C \normeLinf{\omega} \ln (e+ \frac{X^{20} }{\normeLinf{\omega}}),
\end{eqnarray*}
which proves (\ref{eststatique}).

Now let us prove (\ref{estjbord}). As $v  = \restr{\bar{v}}{\Omega} - \nabla\alpha$
(see~page~\pageref{pagequatre}), we have
\begin{eqnarray*}
[ \nabla v ] ^{\Omega'} _r & \leq & \| \nabla \alpha \|_{1+r} + [\nabla \bar{v } ] ^{\Omega'} _r \\
& \leq & C X ^{20} + [\nabla \bar{v } ] ^{\Omega'} _r.
\end{eqnarray*}
Substituting (\ref{stat8}) in (\ref{stat7}), we get
\begin{eqnarray*}
\partial_j\bar{v}^i & = & \Lambda^{-2 } \chi(D) (\partial_j\bar{v}^i + \sum_{k,l =1}^3 \alpha_{k,l }^i a _{jk
} \bar{\omega}^l) \\
& & \mbox{}+ \sum_{k,l =1}^3 \alpha_{k,l } \Lambda^{-2} [\partial_j\partial_{k } \bar{\omega}^l- \Delta (a
_{jk} \bar{\omega}^l)] - \sum_{k,l =1}^3 \alpha_{k,l }a _{jk} \bar{\omega}^l .
\end{eqnarray*}
This directly implies (\ref{estjbord}), because all terms but those in the last summation are bounded in
$C^r(\rtrois)$ by $C X ^{20}$, and
\begin{eqnarray*}
[a _{jk} \bar{\omega}^l]^{\Omega'}_r & \leq & \normeCr{a _{jk}}\normeLinf{\bar{\omega}} + \normeLinf{a
_{jk}}[\omega^l]^{\Omega'}_r \\
& \leq & C X ^{20} + [\omega] ^{\Omega'} _r.
\end{eqnarray*}

\subsection{Dynamic estimates}
\label{dynam}
Now we prove (\ref{dyna13}), (\ref{dyna15}) and (\ref{dyna18}). Let's recall that we want
to control the evolution of $w ^\nu$ in $C^r(\Omega)$, of $\teste{\nabla}{w
^\nu\otimes \omega}$ in $C^{r-1}(\Omega)$ and of $\omega \cdot n$ in $C^r(\bord{\Omega})$.
\subsubsection{Estimation of $w^\nu$ in $C^r(\Omega)$} \label{estimrmi}
We need to estimate $w ^\nu\pgrad v$ in $C^r(\Omega)$. As $\nu\in\{1,\ldots,N'\}$ is fixed,
we shall write, in this section, $w$ for $w^\nu$. The idea is to estimate $ \| w \cdot
\nabla v \|_{r }$ by
the norms of $ \divergence (w \cdot \nabla v)$ and of $\rotationnel (w \cdot \nabla v)$ in
$C^{r-1}(\Omega)$ and the norm of $(w \cdot \nabla v)\cdot n$ in $C^r(\bord{\Omega})$.
\paragraph{Divergence}
We write
\begin{eqnarray*}
\divergence (w \pgrad v) & = & \partial_i(w ^j\partial_jv ^i) \\
& = & \restr{\partial_i(\bar{w} ^j\partial_j\bar{v} ^i)}{\Omega} \\
& = & \restr{\partial_i[\reste{\bar{w }^j}{\partial_j\bar{v }^i}+ \paraproduit{\bar{w
}^j}{\partial_j\bar{v
}^i}+ \paraproduit{\partial_j \bar{v }^i }{\bar{w }^j}]}{\Omega},
\end{eqnarray*}
where $\bar{v}$ et $\bar{w}$ are extensions of $v$ and $w$ respectively, to be specified.
The remainder and the
paraproducts are treated as usual (see~\cite{chemin}, chapter~2):
\[
\| \reste{\bar{w }^j}{\partial_j\bar{v }^i} \|_{r } \leq C \| \bar{w}^j \|_{r } \|
\partial_j\bar{v }^i \|_{0
}, \]
\begin{eqnarray*}
\| \partial_i\paraproduit{\bar{w }^j}{\partial_j\bar{v }^i} \|_{r-1 } & \leq & \|
\paraproduit{
\partial_i\bar{w }^j}{\partial_j\bar{v }^i} \|_{r-1 } + \| \paraproduit{\bar{w
}^j}{\partial_j
\divergence\bar{v }} \|_{r-1 } \\
& \leq & C \normeCrmu{\partial_i\bar{w}^j} \normezer{\partial_j\bar{v}^i} + C
\normeLinf{\bar{w}^j}\normeCrmu{\partial_j\divergence\bar{v}},
\end{eqnarray*}
and
\begin{eqnarray*}
\normeCrmu{\partial_i \paraproduit{\partial_j \bar{v }^i }{\bar{w }^j} } & \leq &
\normeCrmu{\paraproduit{\partial_j \mbox{\scriptsize{\divergence}}\bar{v }}{\bar{w }^j} } +
\normeCrmu{ \paraproduit{\partial_j
\bar{v }^i }{\partial_i\bar{w }^j} } \\
& \leq & C \normeCrmu{\partial_j\divergence\bar{v}}\normezer{\bar{w}^j} + C
\normeLinf{\partial_j\bar{v}^i}
\normeCrmu{\partial_i\bar{w}^j},
\end{eqnarray*}
thus
\[ \normeCrmu{ \divergence( w \pgrad v ) } \leq C \normeCr{\bar{w}}
(\normeLinf{\nabla\bar{v}}+ \normeCrmu{ \nabla
\divergence \bar{v}}).
\]
Choosing for $w$ any continuous extension from $C^r(\Omega)$ into $C^r(\rtrois)$, and
setting $\bar{v}= P_c v$ (see lemma~\ref{lem2}), one obtains
\begin{equation}
\normeCrmu{\divergence(w \pgrad v)} \leq C \normeLip{v } \normeCr{w }.  \label{dyna1}
\end{equation}
\paragraph{Rotational}
Let us consider the third component (for example) of the vorticity:
\begin{eqnarray*}
[\rotationnel(w \pgrad v)]^3 & = & \partial_1 (w ^j\partial_j v ^2) - \partial_2 (w
^j\partial_j v ^1) \\
& = & \restr{ [\partial_1 (\bar{w} ^j\partial_j \bar{v} ^2) - \partial_2 (\bar{w}
^j\partial_j \bar{v} ^1)]
}{\Omega}.
\end{eqnarray*}
Here continuous extensions $w \in C^r(\Omega)$ to $\bar{w}\in C^r(\rtrois)$ and $v
\in\Lip(\Omega)$ to $\bar{v}\in\Lip(\rtrois)$ will be good enough.  It is easily seen
(making use of $T$ and $R$) that
\begin{eqnarray*}
\partial_1 (\bar{w} ^j\partial_j \bar{v} ^2) & = & \partial_1\partial_j(\bar{w}^j \bar{v}^2)
-
\partial_1(\bar{v}^2 \divergence \bar{w}) \\
& = & \partial_j(\bar{w}^j\partial_1\bar{v}^2) + [\partial_j(\bar{v}^2 \partial_1\bar{w}^j)
-
\partial_1(\bar{v}^2 \divergence \bar{w})].
\end{eqnarray*}
Let us consider the difference that has been put between brackets. We estimate separately
each remainder:
\[ \normeCr{ \reste{\partial_1\bar{w}^j}{\bar{v}^2} } \leq C \normeCrmu{\partial_1\bar{w}^j
}
\normeun{\bar{v}^2}, \]
and
\[ \normeCr{ \reste{\bar{v}^2 }{\divergence \bar{w} } } \leq C \normeun{\bar{v}^2}
\normeCrmu{\divergence
\bar{w} }.
\]
When it comes to paraproducts, however, we need to exploit some simplifications:
\begin{eqnarray*}
\lefteqn{ \partial_j\paraproduit{\partial_1\bar{w}^j}{\bar{v}^2} +
\partial_j\paraproduit{\bar{v}^2}{\partial_1\bar{w}^j} -
\partial_1\paraproduit{\bar{v}^2}{\divergence\bar{w}}
- \partial_1\paraproduit{\mbox{\scriptsize{\divergence}}\bar{w}}{\bar{v}^2} } \\
& = & \paraproduit{\partial_1\bar{w}^j}{\partial_j\bar{v}^2} +
\paraproduit{\partial_j\bar{v}^2}{\partial_1\bar{w}^j} -
\paraproduit{\partial_1\bar{v}^2}{\divergence\bar{w}}
- \paraproduit{\mbox{\scriptsize{\divergence}}\bar{w}}{\partial_1\bar{v}^2},
\end{eqnarray*}
and
\[
\normeCrmu{\paraproduit{\partial_1\bar{w}^j}{\partial_j\bar{v}^2}} \leq C
\normeCrmu{\partial_1\bar{w}^j}
\normezer{\partial_j\bar{v}^2},
\]
\[
\normeCrmu{\paraproduit{\partial_j\bar{v}^2}{\partial_1\bar{w}^j}} \leq C
\normeLinf{\partial_j\bar{v}^2}
\normeCrmu{\partial_1\bar{w}^j},
\]
\[
\normeCrmu{\paraproduit{\partial_1\bar{v}^2}{\divergence\bar{w}}} \leq C
\normeLinf{\partial_1\bar{v}^2}
\normeCrmu{\divergence\bar{w}},
\]
\[
\normeCrmu{\paraproduit{\mbox{\scriptsize{\divergence}}\bar{w}}{\partial_1\bar{v}^2}} \leq C
\normeCrmu{ \divergence\bar{w}}
\normezer{\partial_1\bar{v}^2}.
\]
We do the same thing on $\partial_2 (\bar{w}^j\partial_j\bar{v} ^1)$, which gives after
subtraction
\[   [\rotationnel (w \pgrad v)]^3 = \teste{\nabla}{w \otimes \omega^3 } + Y,
\]
with $\normeCrmu{Y}\leq C \normeLip{\bar{v}} \normeCr{\bar{w}}$.  So we get
\begin{equation}
      \normeCrmu{\rotationnel(w \pgrad v)} \leq C (\normeLip{v } \normeCr{w } +
                                                  \normeCrmu{\teste{\nabla}{w \otimes
\omega}}).
                                                  \label{dyna2}
\end{equation}

\paragraph{Normal component}  \label{normcomp}
Since $v$ and $w$ are both tangent to $\bord{\Omega}$ (see p.~\pageref{wnutangents}), one
has
%
%\begin{eqnarray*}
%
$       (w \pgrad v) \cdot n =  - w^j \partial_j v^i \partial_i \delta
%& = & - w^j \partial_j (v^i \partial_i \delta) + w^j v^i \partial_i\partial_j \delta
=   w^j v^i \partial_i\partial_j \delta$:
%\end{eqnarray*}
%
as in~\cite{boubre}, we borrow some smoothness from the boundary, which leads to
\begin{equation}
\normeCr{(w \pgrad v) \cdot n } \leq C (\normeCr{w } \normeLinf{v } + \normeCr{v }
\normeLinf{w }).
\label{dyna3}
\end{equation}
\paragraph{Consequence for the field}
%$w ^\nu$}
The inequalities (\ref{dyna1}), (\ref{dyna2}) and (\ref{dyna3}) imply
\begin{equation}
\normeCr{w ^\nu \pgrad v } \leq C (\normeLip{v } \normeCr{w ^\nu} +
\normeCrmu{\teste{\nabla}{w^\nu \otimes
\omega}}) ;
                                    \label{dyna4}
\end{equation}
the required estimate is proved in section~1.6 of~\cite{mathese}. Combining (\ref{dyna4})
and (\ref{defwnu}), one gets (\ref{dyna13}).
\subsubsection{Estimation of $\left\protect\langle\nabla, \, w ^\nu \otimes \omega
\right\protect\rangle$ in $C^{r-1}(\Omega)$}
Throughout this section $\nu$ is fixed, so we will again drop the $\nu$ in $w ^\nu$.
By  multiplication of (\ref{eqtourb}) and (\ref{defwnu}) (see~\cite{gamstr}, page~402), one
has
\begin{equation}
\left\{\begin{array}{l}
\partial_t(w \otimes \omega) + \teste{\nabla}{v \otimes (w \otimes \omega)} =
        (w \pgrad v) \otimes \omega + w \otimes (\omega \pgrad v), \\
\restr{(w \otimes\omega)}{t=0 } = w _0 \otimes \omegazero.         \label{dyna5}
\end{array}\right.
\end{equation}
Applying $\teste{\nabla}{\cdot}$ to the equation of~(\ref{dyna5}) one gets, $\forall i
\in\{1,\,2,\,3 \}$,
\[
\partial_t\partial_j(w ^j\omega^i)+ \partial_j\partial_k(v ^k w ^j\omega^i) = \partial_j(w
^k \omega^i
\partial_k v ^j) + \partial_j(w ^j \omega^k \partial_k v ^i),
\]
i.e.
\[
\partial_t\partial_j(w ^j\omega^i)+ \partial_k[v ^k \partial_j(w ^j\omega^i)] = \partial_j(w
^j \omega^k
\partial_k v ^i),
\]
so (\ref{dyna5}) implies that $u \egnot \teste{\nabla}{w \otimes \omega} \in \Linf([0,T];
C^{r-1}(\Omega))$
satisfies the system
\begin{equation}
\left\{\begin{array}{l}
\partial_t u + \teste{\nabla}{v \otimes u} =
        \teste{\nabla}{w \otimes (\omega \pgrad v) } , \\
\restr{\teste{\nabla}{w \otimes\omega}}{t=0 } = u _0 \egnot \teste{\nabla}{w _0 \otimes
\omegazero}.
\label{dyna6}
\end{array}\right.
\end{equation}
We have to estimate $f \egnot \teste{\nabla}{w \otimes (\omega \pgrad v) }$ in
$C^{r-1}(\Omega)$. We follow
the proof of~\cite{gamstr}, page 404, changing here and there some details:
\begin{eqnarray}
f ^i & = & \partial_j(w ^j \omega^k \partial_k v ^i) \nonumber\\
& = & \partial_j\partial_k (w ^j \omega^k v ^i) - \partial_j[\partial_k(\omega^k w ^j) v ^i]
\nonumber\\
& = & \partial_k[ \partial_j(w ^j \omega^k) v ^i - \partial_j(\omega^j w ^k) v ^i] +
\partial_k(\omega^k w ^j
\partial_j v ^i). \label{dyna10}
\end{eqnarray}
Let us start with the last term, $\partial_k(\omega^k w ^j \partial_j v ^i)$. We set
$\bar{\omega}= P
\omega$, so that $\normeLinf{\bar{\omega}}+ \normeLinf{\divergence\bar{\omega}} \leq C
\normeLinf{\omega}$, and we extend $w \cdot \nabla v$ continuously from $C^r(\Omega)$ into
$C^r(\rtrois)$. Then we have
\begin{eqnarray*}
\partial_k(\omega^k w ^j \partial_j v ^i) & = & \partial_k \restr{(\bar{\omega}^k
\overline{w \cdot \nabla v }
\,^i)}{\Omega} \\
& = & \partial_k \restr{[ \reste{\bar{\omega}^k}{\overline{w \cdot \nabla v } \,^i} +
\paraproduit{\bar{\omega}^k} {\overline{w \cdot \nabla v } \,^i} + \paraproduit{\overline{w
\cdot \nabla v }
\,^i}{\bar{\omega}^k} ] }{ \Omega },
\end{eqnarray*}
with
\begin{eqnarray}
\normeCr{\reste{\bar{\omega}^k}{\overline{w \cdot \nabla v } \,^i} +
\paraproduit{\bar{\omega}^k} {\overline{w
\cdot \nabla v } \,^i}} & \leq & C \normeLinf{\bar{\omega}^k} \normeCr{\overline{w \cdot
\nabla v } \,^i}
\nonumber\\
& \leq & C \normeLinf{\omega} \normeCr{w \pgrad v }, \label{dyna7}
\end{eqnarray}
and
\begin{eqnarray}
\normeCrmu{\partial_k \paraproduit{\overline{w \cdot \nabla v } \,^i}{\bar{\omega}^k} } &
\leq & \normeCrmu{
\paraproduit{\overline{w \pgrad v }}{\divergence \bar{\omega}} } + \normeCrmu{\paraproduit{
\partial_k
\overline{w \pgrad v }}{\bar{\omega}^k } } \nonumber\\
& \leq & C \normeLinf{\overline{w \pgrad v }} \normeCrmu{\divergence\bar{\omega}} + C
\normeCrmu{\nabla
\overline{w \pgrad v }} \normezer{\bar{\omega}} \nonumber\\
& \leq & C \normeLinf{\omega} \normeCr{w \pgrad v }. \label{dyna8}
\end{eqnarray}
Combining (\ref{dyna7}), (\ref{dyna8}) and (\ref{dyna4}), we find
\begin{eqnarray}
\lefteqn{ \normeCrmu{\partial_k(\omega^k w ^j \partial_j v ^i)} } \nonumber \\
& \leq & C \normeLip{v }\normeLinf{\omega} \normeCr{w } + C \normeLinf{\omega} \normeCrmu{
\teste{\nabla}{w
\otimes \omega} }. \label{dyna9}
\end{eqnarray}
Now consider the first term of~(\ref{dyna10}). We extend continuously $w$ (to $\bar{w}$)
from
$C^r(\Omega)$ into $C^r(\rtrois)$, $v$ (to $\bar{v}$) from $\Lip(\Omega)$ into
$\Lip(\rtrois)$, and
we set $\bar{\omega}= P \omega$. Let $\zeta = \teste{\nabla}{ \bar{w} \otimes \bar{\omega} -
\bar{\omega} \otimes \bar{w} }$. Using these notations, the first term of~(\ref{dyna10}) is
written $
\restr{\teste{\nabla} {\zeta \otimes \bar{v}} ^i}{\Omega}$. Remark that $\divergence \zeta =
0$, so
%\begin{eqnarray*}
$
\teste{\nabla}{\zeta \otimes \bar{v}}
%& = & \partial_k(\zeta^k\bar{v}) \\
%&
 =
%&
\partial_k[\reste{\zeta^k}{\bar{v}} + \paraproduit{\zeta^k}{\bar{v}}] +
\paraproduit{\partial_k\bar{v}}
{\zeta^k},
$
%\end{eqnarray*}
and hence
$
\normeCrmu{\teste{\nabla}{\zeta \otimes \bar{v}}} \leq C \normeLip{\bar{v}}
\normeCrmu{\zeta}.
$
As
\begin{eqnarray*}
%\lefteqn{
\normeCrmu{\zeta}
%}
%\nonumber \\
& \leq &
\normeCrmu{\teste{\nabla}{\bar{w} \otimes \bar{\omega}} } + C \normeLinf{\bar{\omega}}
\normeCr{\bar{w}} + C \normeCrmu{\divergence\bar{\omega}} \normeLinf{\bar{w}}
%\nonumber
\\
& \leq & C ( \normeCrmu{\teste{\nabla}{w \otimes \omega}}+ \normeLinf{w }
\| \omega \cdot n \|_{C^r(\bord{\Omega})} + \normeLinf{\omega} \normeCr{w } ),
%\nonumber\\
%& & \mbox{} %\label{dyna14}
\end{eqnarray*}
which can be shown by similar arguments to those of section~\ref{nablawbomb}, we finally
have on $f$, taking (\ref{dyna9}) into account,
%  and (\ref{dyna14})
the estimate
\begin{eqnarray}
\lefteqn{ \normeCrmu{f } } \nonumber \\
& \leq & C \normeLip{v } (\normeCrmu{\teste{\nabla}{w \otimes \omega}} + \normeLinf{\omega}
\normeCr{w }
\nonumber\\
& & \mbox{}+ \normeLinf{w } \| \omega \cdot n \|_{C^r(\bord{\Omega})}).
\label{dyna11}
\end{eqnarray}
We then extend continuously $f$ to $\bar{f} \in \Linf([0,T]; C^{r-1}(\rtrois))$, $u _0$ to
$\bar{u _0} \in
C^{r-1}(\rtrois)$, $v$ to $\bar{v}$ such that $\normeLip{\bar{v}} +
\normeLip{\divergence\bar{v}} \leq C
\normeLip{v }$ (setting $\bar{v}= P_c v$). Next we approach $\bar{f}$ and $\bar{u_0 }$ by
regularizations with
fixed $t$ ($f _{k} \egnot \conv{\rho_{k^{-1}}}{\bar{f}}$, $u _{0,k} =
\conv{\rho_{k^{-1}}}{\bar{u_0}}$, $k \in
{\mathbb N}$), and we define $u _k$, $\forall k \in {\mathbb N}$, as the solution of
\[
\left\{\begin{array}{l}
\partial_t u _k + \teste{\nabla}{\bar{v} \otimes u _k} = f _k, \\
\restr{u _k}{t=0 } = u _{0,k},
\end{array}\right.
\]
in the space $\Linf([0,T]; \Lip(\rtrois))$. We apply to the $u _k$'s the following lemma,
which is almost identical
to lemma~4.1.1 in~\cite{chemin}.
\begin{lem}
Let $T>0$ and $r \in\,]0,1[$. Suppose that $g$, $h$ and $w$, all three in $\Linf([0,T];
\Lip(\rtrois))$,
satisfy
\[
\partial_t g + \teste{\nabla}{w \otimes g }= h,
\]
and that $\divergence w$ also belongs to $\Linf([0,T]; \Lip(\rtrois))$. Then one has, with
the notation $W (\cdot ) = \normeLip{w (\cdot )} + \normeLip{\divergence w (\cdot )}$, the
inequality
\[
\normeCrmu{g (t)} \leq \normeCrmu{g (0)} e ^{C \int_0 ^t W (s) \,ds} + \int_0^t \normeCrmu{h
(s) }
e ^{C \int _s ^t W (\tau)  \,d\tau},
\]
the constant $C$ depending only on $r$.
\end{lem}
\demonstration $\forall q$ integer $ \geq -1$, one has
\[
\left\{\begin{array}{l}
\partial_t \Delta_q g + w \pgrad \Delta_q g = \deltaq h +  w \pgrad \Delta_q g - \deltaq (w
\pgrad g) -
\deltaq (g \,\divergence w) , \\
\restr{\deltaq g }{t=0 } = \deltaq g (0).
\end{array}
\right.
\]
The inequality
\[
\normeLinf{w \pgrad \Delta_q g - \deltaq (w \pgrad g)} \leq 2^{-q(r-1)} W \normeCrmu{g}
\]
is obtained as in the proof of lemma~4.1.1 in~\cite{chemin}, except that the
commutators $[\deltaq, \reste{\divergence w }{\cdot}]$ (see~page 69 in~\cite{chemin}) are
estimated by taking advantage of the fact that $\divergence w \in \Lip(\rtrois)$:
\begin{eqnarray*}
\lefteqn{ \normeLinf{\reste{\divergence w }{\deltaq g } - \deltaq\reste{\divergence w }{g }}
}\\
& \leq & \normeLinf{\reste{\divergence w }{\deltaq g }} +
\normeLinf{\deltaq\reste{\divergence w }{g }} \\
& \leq &  C \normeLip{\divergence w } \normeLinf{\deltaq g } + C 2^{-q(r-1)}
\normeCrmu{\reste{\divergence w
}{g }} \\
& \leq & C 2 ^{-q(r-1) } \normeLip{\divergence w } \normeCrmu{g }.
\end{eqnarray*}
The additional term does not pose any problems:
\begin{eqnarray*}
\normeLinf{\deltaq (g \,\divergence w) } & \leq & C 2^{-q(r-1)} \normeCrmu{g \,\divergence w
} \\
& \leq & C 2^{-q(r-1) } \normeLip{\divergence w} \normeCrmu{g }.
\end{eqnarray*}
One concludes as in~\cite{chemin} ($\Linf$ estimations of $\deltaq g$, multiplication by
$2^{q(r-1)}$,
supremum over $q$, Gronwall).
\findemonstration

\noindent So we have, $\forall k \in {\mathbb N}$,
\[
\normeCrmu{u _k(t)} \leq \normeCrmu{u _{0,k }} e ^{C \int_0^t \normeLip{v (s)}\,ds} +
\int_0^t \normeCrmu{f
_k(s)} e ^{C\int _s^t\normeLip{v (\tau)}\,d \tau} ds,
\]
and by subtraction, $\forall k,l \in {\mathbb N}$,
\begin{eqnarray*}
\lefteqn{ \normeCrpmu{u _k(t)- u _l(t)} }\\
& \leq & \normeCrpmu{u _{0,k }- u _{0,l }} e ^{C \int_0^t \normeLip{v (s)}\,ds} \\
& & \mbox{}+ \int_0^t \normeCrpmu{f _k(s) - f _{l }(s)} e ^{C\int _s^t\normeLip{v (\tau)}\,d
\tau} ds,
\end{eqnarray*}
with any $r'\in\,]0,r[$. We deduce from this that the sequence $(u_k)_{k \in {\mathbb N}}$
converges in
$\Linf([0,T]; C^{r'-1}(\rtrois))$ to a distribution $\tilde{u} \in\Linf([0,T];
C^{r-1}(\rtrois))$ solution
of
\[
\left\{\begin{array}{l}
\partial_t \tilde{u}  + \teste{\nabla}{\bar{v} \otimes \tilde{u}} = \bar{f}, \\
\restr{\tilde{u}}{t=0} = \bar{u}_0,
\end{array}\right.
\]
with the estimate
\begin{equation}
%
%\lefteqn{
\normeCrmu{\tilde{u} (t)}
%}
%\nonumber \\
%& \leq & \normeCrmu{\bar{u} _0} e ^{C \int_0^t \normeLip{v (s)}\,ds} + \int_0^t
%\normeCrmu{\bar{f} (s)} e
%^{C\int _s^t\normeLip{v (\tau)}\,d \tau}\,ds \nonumber\\
%&
\leq
%&
\normeCrmu{u _0} e ^{C \int_0^t \normeLip{v (s)}\,ds} + \int_0^t \normeCrmu{f (s)}
e ^{C\int
_s^t\normeLip{v (\tau)}\,d \tau} \, ds. \label{dyna12}
\end{equation}
As $\restr{\tilde{u }}{\Omega}$ is a solution of the system (\ref{dyna6}) in $\Linf([0,T];
C^{r-1}(\rtrois))$, $u= \restr{\tilde{u }}{\Omega}$, by uniqueness, so $\normeCrmu{u (t)}$
is not greater than
(\ref{dyna12}), $\forall t \in[0,T]$.
The inequalities (\ref{dyna11}) and (\ref{dyna12}) allow us to conclude:
\begin{eqnarray}
\lefteqn{ \normeCrmu{\teste{\nabla}{w ^{\nu}(t) \otimes \omega(t)}} } \nonumber \\
& \leq & \normeCrmu{\teste{\nabla}{w ^{\nu}_0  \otimes \omega_0}} e ^{C \int_0^t \normeLip{v
(s)}\,ds}
\nonumber\\
& & \mbox{}+C \int_0^t \normeLip{v (s)} (\normeCrmu{\teste{\nabla}{w ^{\nu}(s) \otimes
\omega(s)}} +
\normeLinf{\omega(s)} \normeCr{w ^\nu (s)} \nonumber\\
& & \mbox{} \quad + \normeLinf{w ^\nu (s)} \| \omega(s)\cdot n \|_{C^r(\bord{\Omega})}) e
^{C \int_s^t
\normeLip{v (\tau)}\, d \tau } ds, \nonumber
\end{eqnarray}
thus, dividing by $\normeLinf{\omega(t)}$, we obtain (\ref{dyna15}).
\subsubsection{Estimation of $\omega \cdot n$ in $C^r(\bord{\Omega})$}
\label{evolomn}
It is shown in~\cite{gamstr}, page 406, that $\forall\mu,\nu \in\{1,\ldots,N' \}$,
\begin{equation}
\partial_t(w ^\mu \produitvectoriel w ^\nu) + \teste{\nabla}{v \otimes (w ^\mu
\produitvectoriel w ^\nu)} = -
^t (\nabla\otimes v) (w ^\mu \produitvectoriel w ^\nu), \label{dyna16}
\end{equation}
where $[ ^t (\nabla\otimes v) X ]^j\egnot \sum_{i=1 }^3 X ^i \partial_j v ^i$,
$j=1,\,2,\,3$. It follows
from (\ref{dyna16}) and (\ref{eqtourb}) that $\forall\mu,\nu \in\{1,\ldots,N' \}$, the
quantity $(w ^\mu \produitvectoriel w
^\nu) \cdot \omega$ is preserved along the flow lines:
\begin{eqnarray*}
\lefteqn{ (\partial_t + \teste{\nabla}{v \otimes \cdot}) [(w ^\mu \produitvectoriel w ^\nu)
\cdot \omega] }\\
& = & - [(w ^\mu\produitvectoriel w ^\nu)^i \partial_j v ^i] \,\omega^j + (w
^\mu\produitvectoriel w ^\nu)^i
\omega^j \partial_j v ^i \\
& = & 0.
\end{eqnarray*}
In particular, $\forall x \in\bord{\Omega}$, $\forall t \in[0,T]$, one has
\begin{equation}
[ (w ^\mu \produitvectoriel w ^\nu) (t, \psi(t,x)) ] \cdot \omega(t, \psi(t,x)) = [w
^\mu_0(x)
\produitvectoriel w ^\nu_0(x)] \cdot \omega_0(x). \label{dyna17}
\end{equation}
As $\restr{w ^\mu (t, \cdot)}{\bord{\Omega}} \cdot n =0$, $\forall\mu\in\{1,\ldots,N' \}$,
$\forall t
\in[0,T]$, the vector products $w ^\mu \produitvectoriel w ^\nu$ have the direction of $n$:
\[
(w ^\mu \produitvectoriel w ^\nu) (t, \psi(t,x)) = (-1) ^{\alpha_{\mu,\nu}(t,x)} | (w ^\mu
\produitvectoriel w
^\nu) (t, \psi(t,x)) | \, n(\psi(t,x)),
\]
with $\alpha_{\mu,\nu}(t,x)= 0 \mbox{ or }1$ following the orientation of $w ^\mu
\produitvectoriel w ^\nu$ in
comparison with $n$.
Moreover, (\ref{dyna16}) implies that $\forall t \in[0,T]$, $\forall x \in\bord{\Omega}$, $w
^\mu_0(x)
\produitvectoriel w ^\nu_0(x) =0$ $\Leftrightarrow$ $(w ^\mu \produitvectoriel w ^\nu) (t,
\psi(t,x))=0$.
Thus $\alpha _{\mu,\nu}$ actually depends  only on $x$, so (\ref{dyna17})
can be rewritten
\begin{eqnarray*}
\lefteqn{ (-1)^{\alpha_{\mu,\nu}(x)} |(w ^\mu \produitvectoriel w ^\nu) (t, \psi(t,x)) |
(\omega \cdot n)(t,
\psi(t,x)) }\\
& = & (-1)^{\alpha_{\mu,\nu}(x)} |w ^\mu_0(x) \produitvectoriel w ^\nu_0(x)| (\omega_0\cdot
n) (x),
\end{eqnarray*}
which we obviously simplify by $(-1)^{\alpha_{\mu,\nu}(x)}$. Then we use the fact that $\{w
^\nu; \nu=1,\ldots,N' \}$ is admissible, $\forall t \in[0,T]$ (see~\cite{gamstr},
corollary~4.3, page~406 again): dividing by $ \sum _{\mu< \nu} | w ^\mu \produitvectoriel w
^\nu |$ gives on
$\bord{\Omega}$, $ \forall t \in[0,T]$,
\[
(\omega\cdot n)(t,\psi(t,\cdot)) = \frac{ \sum_{\mu<\nu} | w ^\mu_0 \produitvectoriel w
^\nu_0 |
\,\omegazero\cdot n }{ \sum_{\mu<\nu} |(w ^\mu \produitvectoriel w ^\nu)(t,\psi(t,\cdot))|
},
\]
hence (\ref{dyna18}).

\subsection{Characteristic functions as multipliers
%Multiplication of $C^{r-1}$ distributions, $r \!\in \,]0,1[$, by the characteristic function of
%a bounded Lipschitzian domain
}
\label{multipl}

The whole section is devoted to the proof of the following fact.
\begin{prop} \label{propmultipl}
The pointwise multiplication by the characteristic function of a bounded Lipschitzian domain is
a continuous mapping from $C^{r-1}(\rtrois)$ into $C^{r-1}(\rtrois)$, $\forall r\in\,]0,1[$.
\end{prop}
This property is proved in~\cite{triebel1} (proposition 3.3.2, pages 197 and 198) only for
\cinf\ domains, but a slight modification of earlier arguments in the same book (section 2.8.5,
pages 149 to 153), gives the result with more limited assumptions about the smoothness of the
boundary. Such a refinement is needed in the proof of proposition~\ref{typepoche} and in
section~\ref{nablawbomb}, when we multiply vector fields by the characteristic function of the
patch.

It is here sufficient to show that, $\forall s \!\in\,]0,1[$, the mapping $B^s_{1,2}(\rtrois)
\rightarrow B^s_{1,2} (\rtrois)$ : $f \mapsto \fcaract{\poche} f$ (pointwise multiplication) is
continuous. Indeed, it follows by duality that \fcaract{\poche}\ is a multiplier in
$(B^s_{1,2})'=B^{-s}_{\infty,2}$, $\forall s \in\,]0,1[$ (see~\cite[p.\ 176--180]{triebel1});
then one can conclude by real interpolation (see~\cite[p.\ 64]{triebel1}).

So let \poche{} be a bounded Lipschitzian domain, and let \fcaract{\poche}{} denote its
characteristic function. Given any function $f \in B^s_{1,2}$, look at the decompostion $
\fcaract{\poche} f = \reste{\fcaract{\poche}}{f } + \paraproduit{\fcaract{\poche}}{f } +
\paraproduit{f }{\fcaract{\poche} }$.

As far as the first two terms are concerned, the regularity of \poche{} is irrelevant:
\[
\|\paraproduit{\fcaract{\poche}}{f } \|_{B^s_{1,2}} + \|\reste{\fcaract{\poche}}{f }
\|_{B^s_{1,2}}
\leq C \normeLinf{\fcaract{\poche }} \| f \|_{B^s_{1,2}},
\]
for a constant $C$ depending only on $s$.

But to estimate the second paraproduct, we need to decompose \fcaract{\poche} itself. Since
\poche{} is Lipschitzian, there exists a collection of \cinfc\ functions, $\{\varphi_i;
i=0,\ldots,N\}$, $N \in {\mathbb N}$, with the following properties (see for
example~\cite{stein}, chapter VI, section 3.2):
\begin{itemize}
\item $\forall i\geq 1$, $\poche\cap\support{\varphi_i}$ (or
$\fermeture{\poche}^c\cap\support{\varphi_i}$) is the set of all points lying below the surface
of equation $x_j=\psi_i(x')$, where $\psi_i$ is Lipschitzian, $j=j (i)=$1, 2 or 3 and
$x'=(x_2,x_3)$, $(x_1,x_3)$ or $(x_1,x_2)$ respectively
\item $\sum_{i=0}^N \varphi_i=1$ on \fermeture{\poche}
\item $\support{\varphi_0}\subset\poche$
\item $ \forall i\geq 1$, $\forall x\in\poche\cap\support{\varphi_i}$,
$\distance(x,\bord{\poche})\geq C |\psi_i(x')-x_j|$, with $C=(1+\normeLip{\psi_i }^2)^{-1/2}$.
\end{itemize}
As $\fcaract{\poche}= \varphi_0+\sum_{i=1}^N \varphi_i \fcaract{\poche}$, we just have to
estimate a typical term, $\paraproduit{f }{(\varphi_i\fcaract{\poche})}$. Set $g=\varphi_i
\fcaract{\poche}$. First we write
\begin{eqnarray}
2^{qs} \normeLun{\deltaq\paraproduit{f }{g }} & \approx & 2^{qs} \normeLun{\deltaq (S_{q-1}f
\deltaq g) } \nonumber\\
& \leq & C 2^{qs} \normeLun{ S_{q-1} f \deltaq g } \nonumber\\
& \leq  & C 2^{qs} \sum_{p=-1}^{q-2} \normeLun{ \deltap f \deltaq g }, \label{multipl0}
\end{eqnarray}
then we look more closely at $\normeLun{ \deltap f \deltaq g }$:
\begin{eqnarray}
\lefteqn{ \int _{\rtrois} | \deltap f (x) | \,  | \deltaq g (x) | \,dx } \nonumber\\
& = & \int_{\rdeux} dx'\, \int_{\rr} dx_j \,| \deltap f (x',\,x_j) |  \, | \deltaq g (x',\,x_j)
| \nonumber\\
& \leq & \int_{\rdeux} dx'\, \sup_{x_j} | \deltap f (x',\,x_j) | \sup_{x'} \int_{\rr} dx_j \, |
\deltaq g (x',\,x_j) |.
\label{multipl1}
\end{eqnarray}
Taking into account that the support of the (one-dimensional, $x'$ being fixed) Fourier
transform of $\deltap f (x',x_j)$ is contained in an interval $[-C 2^p,\, C 2^p]$, we obtain
\begin{eqnarray}
\int_{\rdeux} dx'\, \sup_{x_j} | \deltap f (x',\,x_j) | & \leq & \int_{\rdeux} dx'\, C 2^p
\int_\rr dx _j \, | \deltap f (x',\,x_j) | \nonumber \\
& = & C 2^p \normeLun{\deltap f }. \label{multipl2}
\end{eqnarray}
On the other hand, $\deltaq g$ is the convolution product of $g$ by $2^{3q}h(2^q\cdot)$, where
$h \in \schw (\rtrois)$ has a zero integral. So we have, $\forall x' \in\rdeux$,
\begin{eqnarray}
\lefteqn{ \int_{\rr} dx_j \, | \deltaq g (x',\,x_j) | } \nonumber \\
& = & \int_{\rr} dx_j \, | \int_{\rtrois} dy \, h (y) \,[g
(x'-2^{-q}y',\,x_j-2^{-q}y_j)-g(x',\,x_j)] \,| \nonumber\\
& \leq & \int_{\rtrois} dy \, |h (y)| \int_{\rr} dx_j \,|g
(x'-2^{-q}y',\,x_j-2^{-q}y_j)-g(x',\,x_j) |. \label{multipl3}
\end{eqnarray}
Furthermore, when $x'$ and $y = (y',\,y _j)$ are fixed, it is easily seen that
\begin{equation}
\int_{\rr} dx_j \,|g (x'-2^{-q}y',\,x_j-2^{-q}y_j)-g(x',\,x_j) | \leq C 2^{-q}|y|.
\label{multipl4}
\end{equation}
Indeed, the fact that $x$ and $x-2^{-q}y$ belong or don't belong to \poche{} depends only on the
value of $x_j$; set $J_1=\{x_j; x \in\poche\}$ and $J_2=\{x_j; x-2^{-q}y \in\poche\}$. The left
member of (\ref{multipl4}) is the sum of the integrals of $|g(x-2^{-q}y)-g(x)|$ on $J_1^C \cap
J_2$, $J_1\cap J_2^C$ and $J_1 \cap J_2$. Clearly, the integral on $J_1 \cap J_2$ is smaller
than $\normeLip{\varphi_i} \,2^{-q}|y|$ multiplied by the diameter of the support of
$\varphi_i$. Moreover, if $x \in J_1 \cap J_2^C$, $|g(x-2^{-q}y)-g(x)|=\varphi_i (x)$, so only
points $x \in J_1 \cap J_2^C \cap \,\support{\varphi_i}$ contribute to the integral. And for
these points, whose distance to \bord{\poche} is at most $2^{-q}|y|$, we have the inequality
\[
2^{-q}|y| \geq (1+\normeLip{\psi_i }^2)^{-1/2} |\psi_i(x')-x_j|  ;
\]
therefore their $j$th components are inside an interval of length $\leq C 2^{-q}|y|$. Arguing
along the same lines about $J_1^C \cap J_2$ ends the proof of~(\ref{multipl4}).
Substituting successively (\ref{multipl2}), (\ref{multipl3}) and (\ref{multipl4}) in
(\ref{multipl1}), we get
\begin{eqnarray*}
\normeLun{\deltap f \deltaq g } & \leq & C 2^{p-q} \normeLun{\deltap f } \\
& = & C 2^{p(1-s)} 2^{-q} 2^{ps} \normeLun{\deltap f },
\end{eqnarray*}
hence
\begin{eqnarray*}
(\ref{multipl0}) & \leq & C \sum_{p=-1}^{q-2} 2^{ps} \normeLun{\deltap f } 2^{(1-s)(p-q)} \\
& \leq & \frac{C}{(1-2^{s-1})^{1/2} } (\sum_{p=-1}^{q-2} 2^{2ps} \normeLun{\deltap f }^2
2^{(1-s)(p-q)})^{1/2}.
\end{eqnarray*}
So we obtain
\begin{eqnarray*}
\lefteqn{ \sum_{q=1}^{\infty} 2^{2qs} \normeLun{\deltaq\paraproduit{f }{g }}^2 }\\
& \leq & C _s \sum_{p=-1}^\infty 2^{2ps} \normeLun{\deltap f }^2 \sum_{q=p+2}^\infty
2^{(1-s)(p-q)} \\
& = & C_s'\, \| f \|_{B^s_{1,2} },
\end{eqnarray*}
which concludes the proof.

\section{Regularity results}
\label{regular}
First we show, on the assumptions of theorem~\ref{thm1}, that $\omega(t,\cdot )\in
C^s(\psi(t,P)) \cap C^s(\psi(t,\Omega\setminus\fermeture{P }))$. After that we will prove
that the open set $\psi(t,P)$ is of class $C^s$, $\forall t \in[0,T_0]$.
\subsection{Regularity of the curl}
One must evaluate the differences $| \omega(t, x') - \omega(t, y')|$ with $x',\,y' \in
\psi(t,P)$ and with $ x',\,y' \in  \psi(t,\Omega\setminus\fermeture{P})$.

Let $Q$ be any open set whose closure is a subset of $P$ and let $x,\,y \in Q$. Let us
denote by $\psi_k$ the flow of $v _k$, $\forall k \in \naturels$, where the $v _k$'s are
the smooth solutions corresponding to regularized initial data, as in the end of
section~\ref{sketch}, page~\pageref{regularisation}. Given (\ref{eqtourb}), one can write
\begin{eqnarray*}
\lefteqn{ | \omega_k(t, \psi_k(t,x)) - \omega_k(t, \psi_k(t,y))| }\\
& \leq & | \omega_{0k }(t,x) - \omega_{0k }(t,y)| \\
& & \mbox{}+ \int_0^t |(\omega_k\pgrad v _k)(\tau, \psi_k(\tau,x)) - (\omega_k\pgrad v
_k)(\tau, \psi_k(\tau,y))| \,d\tau \, \\
& \leq & \|\omega_{0i }\|_{s } |x-y| ^s + \int_0^t [\omega_k\pgrad v _k(\tau)] ^{Q,k}_s
|\psi_k(\tau,x) - \psi_k(\tau,y)| ^s\,d \tau,
\end{eqnarray*}
at least $\forall k \geq K$, $K$ depending on $Q$, and with the notation
\[
[f(\tau)] ^{Q,k}_s = [f(\tau)]^{\psi_k(\tau,Q)}_s.
\]
The uniform bound on $v _k$, (\ref{theo1}), imply that $\normeLip{\psi^{\pm 1}_k(\tau,
\cdot )} %
%+ \normeLip{\psi_k^{-1 }(\tau, \cdot )}
\leq C e ^{\int_0^{\tau} \normeLip{v _k(\sigma)} \,d \sigma} \leq  C$, $\forall
\tau\in[0,T_0]$,
where the last constant $C$ does not depend on $\tau$, $k$ or $Q$. So we have
\begin{eqnarray}
 [\omega_k(t)]^{Q,k}_s
& \leq & C \|\omega_{0i }\|_{s } + C \int_0^t [\omega_k\pgrad v _k(\tau)] ^{Q,k}_s \,d
\tau \nonumber\\
& \leq & C \|\omega_{0i }\|_{s } + C \int_0^t \normeLinf{\omega_k(\tau)} [\nabla v
_k(\tau)] ^{Q,k}_s \,d \tau \nonumber\\
& & \mbox{}+ C \int_0^t\normeLinf{\nabla v _k(\tau)}[\omega_k(\tau)] ^{Q,k}_s \,d \tau
\nonumber\\
& \leq & C \|\omega_{0i }\|_{s } + C \int_0^t ([\nabla v _k(\tau)] ^{Q,k}_s +
[\omega_k(\tau)] ^{Q,k}_s) \,d \tau \label{theo3}.
\end{eqnarray}
The estimate (\ref{estjbord}) gives $[\nabla v _k(\tau)] ^{Q,k}_s \leq C X _k^{20} +
[\omega_k(\tau)] ^{Q,k}_s$, hence
\[
[\nabla v _k(\tau)] ^{Q,k}_s \leq C + [\omega_k(\tau)] ^{Q,k}_s, \quad
\forall\tau\in[0,t],
\]
thanks to (\ref{theo1}) and (\ref{majxot}).
Using that in (\ref{theo3}), we have
\begin{equation}
[\omega _k (t)] ^{Q,k} _s  \leq  C e ^{Ct }, \quad\forall t \in[0,T_0], \label{theo4}
\end{equation}
by Gronwall's lemma.

Now let $U$ be any open set whose closure is a subset of $Q$. We show that there exists a
subsequence of $(\omega_{k (l)})_{l \in{\mathbb N}}$ which converges to $\omega$,
uniformly on $\Psi(U) \egnot \{(t,\psi(t,U)); t \in [0,T_0]\}$.  If $k$ is great enough
and if $t, \, t' \in [0,T_0]$ are close enough, then $\psi(t,U) \subset \psi_k(t',Q)$.
Indeed, $\forall x \in U$,
\begin{eqnarray*}
\lefteqn{ |\psi(t,x) - \psi_k(t',x)| }\\
& \leq & | \psi(t,x) - \psi_k(t,x)| + | \psi_k(t,x) - \psi_k(t',x)| \\
& \leq & (\int_0^t \normeLinf{v(\tau,\cdot )-v _k(\tau,\cdot )}\,d\tau) \, e
^{\int_0^t\normeLip{v (\tau,\cdot )}\,d\tau} + | \int_t^{t'} \normeLinf{v _k(\tau)}
\,d\tau | \\
& \leq & C \| v - v _k \|_{\Linf([0,T_0]; \Linf(\Omega))} + C|t-t'|,
\end{eqnarray*}
while the distance between $\psi_k(\cdot ,\,U)$ and $\Omega\setminus\psi_k(\cdot ,\,Q) =
\psi_k(\cdot ,\,\Omega\setminus Q)$ is, on $[0,T_0]$, greater than a positive constant
independent of $k$. Consider, $\forall y,\, y' \in U$, the differences
\begin{eqnarray}
\lefteqn{ | \omega _k(t,\psi(t,y))- \omega_k(t', \psi(t',y')) |} \nonumber \\
& \leq & | \omega _k(t,\psi(t,y))- \omega_k(t', \psi(t,y)) | \label{ettheo1} \\
& & \mbox{} + | \omega _k(t',\psi(t,y))- \omega_k(t', \psi(t',y'))|.
\label{ettheo2}
\end{eqnarray}
One may suppose that there exists a point $x \in Q$ such that $\psi_k(t,x)= \psi(t,y)$;
then (\ref{ettheo1}) is equal to
\begin{eqnarray}
\lefteqn{ | \omega_k(t, \psi_k(t,x)) - \omega_k(t',\psi_k(t,x)) | } \nonumber\\
& \leq & | \omega_k(t, \psi_k(t,x)) - \omega_k(t',\psi_k(t',x)) |
+ | \omega_k(t', \psi_k(t',x)) - \omega_k(t',\psi_k(t,x)) | \nonumber\\
& \leq & | \int_t^{t'} \normeLinf{\omega_k(\tau)} \normeLinf{\nabla v _k (\tau)} \,d \tau|
+ [\omega_k(t')]^{Q,k}_s |\psi_k(t',x) - \psi_k(t,x)| \nonumber\\
& \leq & C|t-t'|, \label{theoun}
\end{eqnarray}
thanks to (\ref{theo4}). Likewise, as soon as $\psi(t,y),\, \psi(t',y') \in \psi_k(t',Q)$,
(\ref{theo4}) gives
\begin{eqnarray}
(\ref{ettheo2}) & \leq & C | \psi(t,y) - \psi(t',y')|^s \nonumber\\
& \leq & C \normeLip{\psi}^s (|t-t'| + |y-y'|) ^s.                      \label{theodeux}
\end{eqnarray}
One deduces from (\ref{theoun}) and (\ref{theodeux}), by Ascoli-Arzela's theorem, the
existence of a subsequence of $(\omega_{k (l)})_{l \in{\mathbb N}}$ that converges to
$\omega$, uniformly on $\Psi(U)$; as $U \fermeturedans Q$ and $Q \fermeturedans P$ are
arbitrary, this implies that $\omega\in C^s(\psi(t,P))$, $\forall t \in[0,T_0]$, because
\begin{eqnarray*}
\lefteqn{ | \omega(t,x) - \omega(t,y)| }\\
& \leq & | \omega(t,x) - \omega_k(t,x)| + | \omega_k(t,x) - \omega_k(t,y)| + |
\omega_k(t,y) - \omega(t,y)| \\
& \leq & 3 C e ^{CT_0 } |x-y| ^s, \quad \forall x, y \in \psi(t,P),
\end{eqnarray*}
for some suitable $k = k(x,y)$. One shows that $\omega \in
C^s(\psi(t,\Omega\setminus\fermeture{P}))$ in the very same way.
\subsection{Regularity of the patch}
We have, $\forall k= k (l)$, $l \in\naturels$,
\[
\left\{\begin{array}{l}
\partial_t w ^\nu_k + \teste{\nabla}{v _k\otimes w ^\nu_k} = w ^\nu_k\pgrad v _k, \\
w ^\nu_k(0) = w ^\nu_0,
\end{array}\right.\quad \nu=1,\ldots,N',
\]
with $\normeLip{v _k(t,\cdot )} + \normeCs{w ^\nu_k(t,\cdot )}$ uniformly bounded both on
$[0,T_0]$ and with respect to $k$. Thus the $w ^\nu_k$'s form, $\forall\nu\in\{1,\ldots,N'
\}$ an equicontinuous sequence on $[0,T_0]\produitcartesien\Omega$: the equicontinuity
with $t$ fixed is clear, and, $\forall k \in \naturels$, $\forall x \in\Omega$, $\forall
t,\,t' \in[0,T_0]$,
\begin{eqnarray*}
\lefteqn{ |w ^\nu_k(t,\psi_k(t,x)) - w ^\nu_k(t',\psi_k(t,x))| }\\
& \leq & |w ^\nu_k(t,\psi_k(t,x)) - w ^\nu_k(t',\psi_k(t',x))|
+ |w ^\nu_k(t',\psi_k(t',x)) - w ^\nu_k(t',\psi_k(t,x))| \\
& \leq & \normeLinf{w ^\nu_k}\normeLinf{\nabla v _k}|t-t'|+ \normeCs{w ^\nu_k(t',\cdot
)}\normeLinf{v _k}^s |t-t'| ^s.
\end{eqnarray*}
Extracting another subsequence, we get, $\forall \nu\in\{1,\ldots,N'\}$, a field $w ^\nu
\in\Linf([0,T_0]; C^s(\Omega))$ solution of
\[
\left\{\begin{array}{l}
\partial_t w ^\nu + \teste{\nabla}{v \otimes w ^\nu} = w ^\nu\pgrad v, \\
w ^\nu(0) = w ^\nu_0.
\end{array}\right.
\]
Finally we deduce the regularity of $\psi(t,P)$ as in~\cite{gamstr}, pages~417 and 418
(theorem~6.1): let $f \in \cunplusr(\rtrois)$ be a function such that $\restr{f }{\bord{P
}}=0 $ and $\nabla f(x) \neq 0$, $\forall x \in \bord{P }$; an equation of $\psi(t,P)$ is
$\varphi(t,x)=0$, where $\varphi \in \Linf([0,T_0]; \Lip)$ is a solution of
\[
\left\{\begin{array}{l}
\partial_t\varphi + \teste{\nabla}{v \otimes\varphi} = 0, \\
\restr{\varphi}{t=0 } = f,
\end{array}\right.
\]
and it is enough to prove that $\nabla\varphi(t,\cdot ) \in C^s$ in a neighbourhood of
$\psi(t,\bord{P})$, $\forall t \in [0,T_0]$. We have
\[
\left\{\begin{array}{l}
\partial_t\nabla\varphi + \teste{\nabla}{v \otimes\nabla\varphi} = - ^t(\nabla\otimes
v)\nabla\varphi, \\
\restr{\nabla\varphi}{t=0 } = \nabla f,
\end{array}\right.
\]
and, $\forall \mu, \,\nu \in\{1,\ldots,N'\}$,
\[
\left\{\begin{array}{l}
\partial_t (w ^\mu \produitvectoriel w ^\nu) + \teste{\nabla}{v \otimes (w
^\mu\produitvectoriel w ^\nu)} = - ^t(\nabla\otimes v) (w ^\mu\produitvectoriel w ^\nu),
\\
\restr{(w ^\mu \produitvectoriel w ^\nu)}{t=0 } = w ^\mu_0 \produitvectoriel w ^\nu_0.
\end{array}\right.
\]
Wherever $\nabla f \neq 0$, one can write
\[
w ^\mu_0 \produitvectoriel w ^\nu_0 = \frac{|w ^\mu_0 \produitvectoriel w ^\nu_0|}{ |
\nabla f |} (-1) ^{\alpha_{\mu,\nu}} \nabla f,
\]
with $\alpha_{\mu,\nu} = 0$ or $1$, depending on the respective orientations of $w ^\mu_0
\produitvectoriel w ^\nu_0$ and $\nabla f$. So in a neighbourhood $V (t)$ of
$\psi(t,\bord{P})$,
\[
(w ^\mu\produitvectoriel w ^\nu)(t,\cdot ) = \frac{|(w ^\mu_0 \produitvectoriel w
^\nu_0)(\psi^{-1 }(t,\cdot ))| }{|(\nabla f)(\psi^{-1 }(t,\cdot ))|} (-1)
^{\alpha_{\mu,\nu}} \nabla\varphi,
\]
hence, by addition,
\[
\nabla\varphi(t,\cdot ) = \frac{|(\nabla f)(\psi^{-1 }(t,\cdot ))|}{\sum_{\mu<\nu}|(w
^\mu_0 \produitvectoriel w ^\nu_0)(\psi^{-1 }(t,\cdot ))|} \sum_{\mu<\nu}(-1)
^{\alpha_{\mu,\nu}}(w ^\mu\produitvectoriel w ^\nu)(t,\cdot ),
\]
so $\nabla\varphi(t,\cdot )\in C ^s(V (t))$, because $\{w ^\nu_0; \nu=1,\ldots,N'\}$ is
admissible and $w ^\nu(t,\cdot )$ $\in$ $C ^s(V (t))$, $\forall \nu \in\{1,\ldots,N'\}$.

\section{Application to 2-D and axisymmetric flows}
\label{twodim}
In 2-D and axisymmetric flows, the vorticity is always tangent to the boundary of the domain. Therefore,
the sum of (\ref{dyna13}) and (\ref{dyna15}) gives
\begin{align*}
X (t) & \leq  C_0 e ^{C\int_0^t \normeLip{v (s)} } \,ds \\
&  + C_0 \int_0^t \normeLip{v (s)}  X (s) \, e ^{C \int_s^t \normeLip{v (\tau)} \,d\tau} \,ds.
\end{align*}
so this time we have, instead of (\ref{Deltath}),
\begin{equation}
\normeLip{v(t)} \leq C (1+\normeLinf{\omega(t)}) \ln(e+ C_0 e ^{C_0 \int_0^t \normeLip{v (\tau, \cdot )}
\,d \tau}).    \label{eststat2dax}
\end{equation}

In the axisymmetric case, as $\normeLinf{\omega(t)}\leq \normeLinf{\frac{\omega(t)}{\delta}} \max_{x \in
\domec} \delta(x)$ $=$ $C \normeLinf{\frac{\omega_0}{\delta}}$ $\max_{x \in \domec} \delta(x),$
(\ref{eststat2dax}) leads to $\normeLip{v (t)} \leq C_0 e ^{C_0 t }$, $\forall t \in\rr^+$, so $T_0$ in
proposition~\ref{propregtang} is arbitrary.

Of course the 2-D case is even more simple, since $\normeLinf{\omega(t)}=\normeLinf{\omega_0}$.

\vspace{1\baselineskip}
\begin{center}
\textbf{REFERENCES}
\vspace{-1\baselineskip}
\end{center}

\footnotetext{This work was supported by a F.R.I.A.\ Grant (Fonds pour
  la Recherche dans l'Industrie et l'Agriculture) when the author was
  a PhD student at the University of Brussels, in Belgium. Last
  revision in November 2001.}
\end{document}